\documentclass[11pt,twoside]{article}
\usepackage{graphicx}
\usepackage{soul}
\usepackage[T1]{fontenc} 
\usepackage{helvet}
\usepackage[margin=0.9in, top=30mm]{geometry}
\usepackage{textcomp}
\usepackage{gensymb}
\usepackage{amsmath}
\usepackage{algorithm,algorithmic}
\usepackage{psfrag}
\usepackage{amsfonts}
\usepackage{amssymb}
\usepackage{color}
\usepackage{color}
\usepackage{graphicx,subfigure}
\usepackage{array}
\usepackage{tabularx}
\usepackage{longtable}
\usepackage{ulem}
\usepackage{cancel}
\usepackage{url}
\definecolor{airforceblue}{rgb}{0.36, 0.54, 0.66}

\newcommand{\revP}[1]{\textcolor{black}{#1}}

\usepackage{booktabs} 

\usepackage{lettrine} 

\usepackage{enumitem} 
\setlist[itemize]{noitemsep} 

\usepackage{abstract} 

\usepackage{fancyhdr} 
\usepackage{titling} 

\pretitle{\begin{center}\huge} 
\posttitle{\end{center}} 
\title{Implementation and Evaluation of\\  Breaking Detection Criteria for a\\ Hybrid Boussinesq  Model}

\author{\textsc{Paola Bacigaluppi} \thanks{
              Institute of Mathematics - University of Zurich, Winterthurerstrasse 190, 8057 Z\"urich, Switzerland, 
        paola.bacigaluppi@math.uzh.ch }
           \and
          \textsc{Mario Ricchiuto} \thanks{Team CARDAMOM, Inria Bordeaux Sud-Ouest, 200 avenue  de la Vieille Tour, 33405 Talence Cedex, France,   mario.ricchiuto@inria.fr  }
                \and \textsc{ Philippe Bonneton} \thanks{ University of Bordeaux, CNRS, UMR 5805 EPOC,  All\'ee Geoffroy Saint-Hilaire, 33615 Pessac, France, philippe.bonneton@u-bordeaux.fr}}

\begin{document}

\maketitle

\begin{abstract}
The aim of the present work is to develop a model able to represent the propagation and transformation of waves in nearshore areas. 
The focus is on the phenomena of wave breaking, shoaling and run-up.
These different  phenomena are represented through a hybrid approach obtained by the coupling of non-linear Shallow Water equations with the extended Boussinesq equations of Madsen and S{\o}rensen.
The novelty is the switch tool between the two modelling equations: a critical free surface Froude criterion. This is based on a physically meaningful new approach to detect wave breaking, 
which corresponds to the steepening of the wave's crest which turns into a roller. 
To allow for an appropriate discretization of both types of equations, we consider a finite element Upwind Petrov Galerkin method with a novel limiting strategy, that guarantees the 
preservation of smooth waves as well as the monotonicity of the results in presence of discontinuities. \\
We provide a detailed discussion of the implementation of the newly proposed detection method, as well as of two other well known criteria which are used for comparison. 
An extensive benchmarking on several problems involving different wave  phenomena and breaking conditions allows to show the robustness of the numerical method proposed, as well as to assess the 
advantages  and limitations of the different detection methods.

\end{abstract}

\label{intro}
\section{Introduction}

Wave breaking is a fundamental  phenomenon, especially in the  
nearshore region. It   drives  wave set-up and set-down,   wave run-up,  nearshore circulations (e.g. longshore and rip currents), sediment transport, forces on structures, and many other effects.
At large scales, an important effect of   wave breaking is the energy dissipation, which is the result of many local hydrodynamic phenomena, ultimately leading to dissipation due to small scale viscous effects.  
In this paper we focus on the wave breaking closure  when wave propagation is performed by means of   depth averaged Boussinesq models.


Among the main families of closures  found in literature we mention the eddy viscosity formulations \cite{zelt1991run,Schaffer1993,Kennedy2000,Sorensen2004,Cienfuegos2010},
the breaking roller methods \cite{BMBBF04,vmf15}, and the hyperbolic or hybrid  approaches \cite{Bonneton2007,Tonelli2011,Tissier2012,Kazolea14}.
Within the first type of approaches,  the eddy viscosity one,   the macroscopic dissipation of the energy is modelled   by means of explicitly  introducing a  viscous dissipation term in the momentum, or, alternatively, both in the mass and momentum  equations.
In particular, the eddy viscosity coefficients are defined via physical arguments, or through some auxiliary evolution model   (see e.g. \cite{kazolea2018} and references therein),
and, eventually, adjusted by means of numerical experiments. An alternative to the eddy viscosity type modeling are the roller approaches, that account more explicitly for the large scale  effects of vorticity on the mean flow.
Finally, the hybrid approaches exploit the properties of hyperbolic conservation laws endowed with an entropy, and model large scale effects of wave breaking with
the dissipation of the total energy across shocks arising in shallow water simulations \revP{or augment/modify Boussinesq-type equations to include, for example, additional  terms caused  by  the presence  of the  surface  rollers (see cf. \cite{SKOTNER1999905})}.\\
Independently of the chosen closure  strategy,   some breaking detection criteria is very often required to trigger
the onset of wave breaking.  In practice it is this criteria that leads to the activation of the closure. 
As well explained in \cite{Okamoto2006}, historically the first detection criteria were based on quantities
computed using information over one full phase of the wave. 
The possibility of performing accurate phase resolved simulations, has led to the necessity of formulating
new detection criteria based on local flow features   \cite{Schaffer1993,Kennedy2000,Sorensen2004,Kazolea14}. 
More specifically, these  criteria rely on a wave-by-wave analysis  more efficient to program in the context of phase resolved simulations, and 
providing a physically  correct detection of breaking onset and termination, provided that   a sufficiently accurate description of the  wave profiles is available.
In this paper we consider a one-to-one comparison of three different detection criteria, tested on a common solver. The wave propagation is handled by the well known  enhanced Boussinesq equations of Madsen and S{\o}rensen \cite{Madsen1992}. 
Wave breaking is modelled using a hybrid approach and reverting locally to the non-linear Shallow Water equations. 
Our focus is devoted to the comparison of numerical implementations of the non-linearity criterion by \cite{Tonelli2009,Tonelli2010,Tonelli2011},
of  the \revP{hybrid slope/vertical velocity criterion by \cite{Kazolea14,ROEBER20121}}, and  of a convective criterion,  which compares the free surface velocity to the celerity of the wave \cite{Bjorkavag2011}. 
 We propose in particular an implementation in an upwind Finite Element solver, with a novel discontinuity capturing technique, ad-hoc developed for this application.
We perform a very thorough benchmarking on flows involving different  types of waves and of breaking configurations.
The results allow to validate our numerical implementation, and    insight  on the advantages and limitations of the physical behaviour of the different detection criteria.
This paper extends and improves the discussion  of  \cite{bacigaluppi:hal-01087945,bacigaluppi:hal-00990002}.\\
The manuscript is organized as follows. In section \S2 we discuss the model equations, while section \S3 is devoted to the discussion of 
 the finite element scheme used to solve them. Some details 
on the  detection and capturing of shocks are provided, as they have been developed ad-hoc for this application, and  provide an  element of novelty. 
Section \S4 encompasses the definition  of the breaking detection criteria, and provides a precise description of their numerical implementation.
Finally, in section \S5 we collect a thorough comparison of the detection criteria on different kinds of waves, along to experimental data.
The paper is concluded by a discussion and outlook on future work. 
\section{Mathematical Model}
\label{model}
\subsection{Enhanced Boussinesq  Equations}
The underlying enhanced Boussinesq model used in this work is the well known system by Madsen and S{\o}rensen reading \cite{Madsen1991}
\begin{equation}
\label{eq:Bsqsys_eq}
\begin{split}
  &\begin{cases}\partial_{t} \eta  + \partial_{x} q = 0\\[5pt]
 \partial_{t} q  -  \mathcal{D}(\eta,q) + \partial_{x} (uq  +  g H^2/2) + gH \partial_{x} h +  C_f q  = 0 
\end{cases} \\
&\text{with} \;\; \mathcal{D}(\eta,q) =   B h^2 \partial_{xxt} q + \beta g h^3 \partial_{xxx} \eta + h \partial_x h \partial_{xt} q/3 + 2 \beta g h^2 \partial_x h \partial_{xx} \eta 
\end{split}
\end{equation}
Here, following Figure \ref{fig_nome},  
$ h$ denotes the depth at still water level,   $\eta$ the free surface level,  $u$ the depth-averaged velocity,  $H=\eta - h\ge0$ the water depth, 
$q=Hu$ the volume flux, and   $g$ the gravitational acceleration.
\begin{figure}[h]
\centering
\includegraphics[trim={0 1.8cm 0 0.4cm},width=0.55\textwidth]{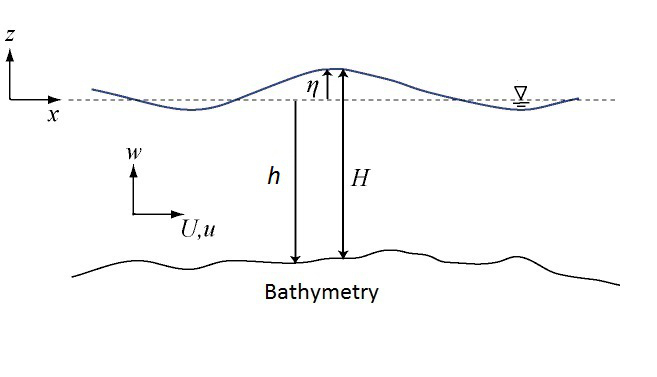}
\caption{Sketch of the main notation.}\label{fig_nome}
\end{figure}
 The coefficient   $C_f$ controls the friction, and  
all dispersive effects are accounted for by the term $\mathcal{D}(\eta,q)$. In particular, the
coefficient values are taken as $\beta=1/15$ and $B=\beta+1/3$. This choice allows to match (in absence of friction) a fourth order Pad\'e development of the dispersion relation of the  Euler equations \cite{Madsen1991}.
\subsection{Wave Breaking Closure: a Hybrid Boussinesq Type Model }
To account also for the effects of wave breaking, system \eqref{eq:Bsqsys_eq} must be coupled to some closure model. 
In the following, we consider a hybrid approach consisting in reverting locally to the non-linear Shallow Water equations, obtained from \eqref{eq:Bsqsys_eq}  with $\mathcal{D}=0$.
The rationale behind this approach is that  in the region of shallow waters, the waves  will locally steepen and rapidly turn into  shocks, across which the volume flux and the momentum are conserved, while the
total energy is dissipated.  The interested reader can refer to  \cite{Bonneton2007,Tonelli2009,Kazolea14,Tissier2012,kazolea2018} for further justification. 

To implement this closure, we modify the model equations \eqref{eq:Bsqsys_eq} to 
 \begin{equation} \label{eqC}
 \begin{cases}
 \partial_{t} \eta  + \partial_{x} q = 0\\[5pt]
\partial_{t} q  - f_{\text{break}}(x,t)  \mathcal{D}(\eta,q) + \partial_{x} (uq  +  g H^2/2) + gH \partial_{x} h + C_f q  = 0, \\
  \end{cases}
   \end{equation}
with the definition of $\mathcal{D}(\eta,q)$ as shown in \eqref{eq:Bsqsys_eq}. The newly introduced term $f_{\text{break}}$, denotes a so-called flag, that  allows to switch from a propagation region ($f_{\text{break}}=1$) to a breaking one ($f_{\text{break}}=0$). \\
For planar waves (one spatial dimension) there are two key elements to define $f_{\text{break}}$:  
the detection of breaking onset and termination; the definition of the region in which the breaking closure is activated.
This paper focuses mainly on the study of the detection criteria. Note that, although the analysis is performed in the context of hybrid wave breaking closure,
detection criteria are an essential element for other approaches as well.
The reader might refer to the extensive references provided in \cite{Kazolea14,kazolea2018} for further details.
Our results will also be relevant for multi-dimensional flows, which will however require to embed some
additional elements related to the direction of wave propagation, and on  water level variations in the transversal direction.


\section{Numerical Approximation}
\label{sec:numericalmeth}
\subsection{Non-linear Upwind Finite Element Discretization}
\label{completeEq}
We propose a non-linear  variant of the upwind stabilized finite element method discussed in \cite{rf14}.
The method exploits a standard   $C^0$ finite element approximation on  a tessellation of the spatial domain $\Omega=[0, L]$,
consisting of $N$ constant intervals, or cells. The cell width is denoted by $\Delta x$.
Let us introduce the vector of unknowns  $ \mathbf{W}=[ \eta\, ,    q ]^T$, along with
 $$
 F(\mathbf{W})=\left[\begin{array}{c} q \\ uq+ g \dfrac{H^2}{2}\end{array}\right]\!,\;
 S_h=gH\left[\begin{array}{c}0\\ \partial_x h \end{array}\right]\!,\;
  S_f=\left[\begin{array}{c}0\\ C_f q \end{array}\right]\!,\;
 D=\left[\begin{array}{c}0\\\mathcal{D}(\mathbf{W})\end{array}\right].
 $$
Denoting by $\varphi_i(x)$ (i=0, ...,N) the standard  linear Lagrange  basis functions, the method relies on  the expansions  
\begin{align}
\mathbf{W}_{\textrm{h}}(t,x)  =\sum_{j=0}^{N} \varphi_j(x)\mathbf{W}_j(t)\,,\quad 
F_{\textrm{h}}(t,x)=  \sum_{j=0}^{N} \varphi_j(x) F(\mathbf{W}_j(t)) .
\label{W_discret}
\end{align}
Introducing the residual
$$
\textrm{r}_{\textrm{h}}=\partial_t\mathbf{W}_{\textrm{h}}  -  f_{\text{break}} D(\mathbf{W}_{\textrm{h}} )+\partial_xF_{\textrm{h}}+S_h(\mathbf{W}_{\textrm{h}} )+S_f(\mathbf{W}_{\textrm{h}} )\,,
$$
The discrete equations  at a node $i$ are obtained from the following stabilized variational statement
\begin{equation}
\begin{split}
\int\limits_{\Omega}\varphi_i\textrm{r}_{\textrm{h}}\,dx + \sum_{i=0}^{N-1}\int\limits_{x_i}^{x_i+\Delta x}\left(A(\mathbf{W}_{\textrm{h}})\partial_x\varphi_i\,\tau_{\textrm{s}}\right)\textrm{r}_{\textrm{h}}\,dx =0,
\end{split}
\end{equation}
where $A=\partial F/\partial \mathbf{W}$ denotes the shallow water flux Jacobian, and $\tau_{\textrm{s}}$ a stabilization matrix.  For the latter it is well known (cf. \cite{rf14} and references therein for details), that
by defining the matrix  as\footnote{Matrix absolute values are computed, as usual, by means of an eigen-decomposition.}  $\tau_{\textrm{s}}=|A|^{-1}\Delta x/2$, the resulting scheme can be recast
as the following upwind splitting method with a non-diagonal mass matrix 
\begin{equation}\label{eq:scheme}
\begin{split}
 \sum\limits_{l=-1,0,1}&M_{i,i+l} \Big(\dfrac{d\mathbf{W}_{i+l}}{dt} + S_{f_{i+l}}\Big) +\dfrac{\mathrm{I}-\mathrm{sign}\big(A(\mathbf{W}_{i+1/2})\big)}{2}\phi_{i+1/2}\\
&\qquad\qquad\qquad\qquad\qquad \;\;\;+\dfrac{\mathrm{I}+\mathrm{sign}\big(A(\mathbf{W}_{i-1/2})\big)}{2}\phi_{i-1/2} =\\
 f_{\text{break}_i}&\int\limits_{\Omega}\varphi_iD(\mathbf{W}_{\textrm{h}}) dx+ \sum_{i=0}^{N-1}  f_{\text{break}_{i+1/2}}\int\limits_{x_i}^{x_{i+1}}\left(A(\mathbf{W}_{\textrm{h}})\partial_x\varphi_i\,\tau_{\textrm{s}}\right)  D(\mathbf{W}_{\textrm{h}})dx,
\end{split}
\end{equation}
where we have introduced  the hydrostatic fluctuations $\phi$, with $\phi_{i+1/2}=F_{i+1}-F_{i}+gH_{i+1/2}[0,\, h_{i+1}-h_i]$ and $\phi_{i-1/2}$ approximated analogously. 
The terms  $\mathbf{W}_{i\pm1/2}$ denote  cell arithmetic averages.

As for the right hand side terms,  neglecting the boundary conditions, which will be discussed later, these are evaluated using the following approximations 
(cf. definition of $D$ in the beginning of the section and of $\mathcal{D}$ in \eqref{eq:Bsqsys_eq})
\begin{equation*}
\begin{split}
\int\limits_{\Omega}\varphi_i\mathcal{D}(\mathbf{W}_{\textrm{h}}) dx =
\int\limits_{x_{i-1}}^{x_{i+1}}\big(-&\partial_x(B h^2_{\textrm{h}}\varphi_i) \partial_{xt}q_{\textrm{h}} + \varphi_i \dfrac{1}{3}h_{\textrm{h}}\partial_xh_{\textrm{h}}\partial_{xt}q_{\textrm{h}}\\
+&\beta g   \varphi_i(  h_{\textrm{h}}^3 \partial_{x}w^{\eta}_{\textrm{h}} + 2  h_{\textrm{h}}^2 \partial_xh_{\textrm{h}}\,w^{\eta}_{\textrm{h}})
\big)dx \,,
\end{split}
\end{equation*}
and
\begin{equation*}
\begin{split}
\int\limits_{x_i}^{x_{i+1}}&\left(A(\mathbf{W}_{\textrm{h}})\partial_x\varphi_i\,\tau_{\textrm{s}}\right)  D(\mathbf{W}_{\textrm{h}})dx
 = -  \frac{\text{sign}\big(A(\mathbf{W}_{i+\frac{1}{2}})\big)}{2} D_{i+1/2}, \quad \text{with}\\
 D_{i+1/2}= \int\limits_{x_i}^{x_{i+1}}&\left[
 \begin{array}{c}
 0\\[10pt]
 Bh_{\textrm{h}}^2\partial_{xt}w^q_{\textrm{h}}+\beta g h_{\textrm{h}}^3\partial_{x}w^{\eta}_{\textrm{h}}+h_{\textrm{h}}\partial_xh_{\textrm{h}}\partial_{xt}q_{\textrm{h}}/3+2\beta g h^2_{\textrm{h}}w^{\eta}_{\textrm{h}}
 \end{array}
 \right]dx.
\end{split}
\end{equation*}
The above expressions are readily evaluated by numerical quadrature in each cell using the linear finite element expansions to interpolate all the quantities
with the ``$\textrm{h}$'' subscript. Note that, for any quantity $u$ we have  $\partial_xu_{\textrm{h}}= u_{i+1}-u_i/\Delta x $ over a cell $[x_i,\, x_{i+1}]$.
The auxiliary variables $w^{\eta}$  and $w^q$ are introduced to handle the higher order derivatives. Their values are obtained by locally inverting the projections
	\begin{equation}
\int\limits_{\Omega}\varphi_i w^{\eta}_{\textrm{h}} = - \int\limits_{i-1}^{i+1} \partial_x\eta_{\textrm{h}}\partial_x\varphi_i dx,\,\;\;\;\;
\int\limits_{\Omega}\varphi_i w^{q}_{\textrm{h}}  = - \int\limits_{i-1}^{i+1} q_{\textrm{h}}\partial_x\varphi_i dx.
\end{equation}
 The interested reader can refer to \cite{rf14} or more details. 
%
\subsection{Shock Capturing, Limiters, and Entropy Fix}
In absence of friction and of dispersive effects, if the mass matrices on the left hand side of \eqref{eq:scheme} are lumped, we end up with the first-order upwind flux splitting approximation
\begin{equation*}
\begin{split}
\Delta x \frac{d\mathbf{W}_i}{dt}  + \dfrac{\mathrm{I}-\mathrm{sign}\big(A_{i+1/2}\big)}{2}\phi_{i+1/2}+\dfrac{\mathrm{I}+\mathrm{sign}\big(A_{i-1/2}\big)}{2}\phi_{i-1/2} =0
\end{split}
\end{equation*}
with the short notation $A_{i\pm1/2}=A(\mathbf{W}_{i\pm1/2})$.
The first-order scheme above is the basis of many discontinuity capturing approximations for shallow water flows (see e.g. \cite{Castro2005,Delis2008,CEA20123317}).
The strategy proposed here is to introduce a cell-wise mass-matrix limiter, that allows to switch from this first  order  non-oscillatory method across shocks,
 to the full finite element approximation,  which has been shown to be third-order accurate and is appropriate for propagation \cite{rf14}. 
 Taking into account the explicit form of the mass matrix entries (see  \cite{rf14} for details),   the full non-linear approximation  can 
be written as 
\begin{equation*}
\begin{split}
\Delta x \frac{d\mathbf{W}_i}{dt}   &
+ \frac{\Delta x\,\delta_{i-\frac{1}{2}}}{2}\Big\{\frac{1}{3} [ \frac{d\mathbf{W}_{i-1}}{dt} -  \frac{d\mathbf{W}_{i}}{dt}] + \text{sign}(A_{i-\frac{1}{2}}) \dfrac{d\mathbf{W}_{i-\frac{1}{2}} }{dt}\Big\} \\ 
&+ \frac{\Delta x\,\delta_{i+\frac{1}{2}}}{2}\Big\{\frac{1}{3} [\frac{d\mathbf{W}_{i+1}}{dt} -  \frac{d\mathbf{W}_{i}}{dt}] - \text{sign}(A_{i+\frac{1}{2}}) \dfrac{d\mathbf{W}_{i+\frac{1}{2}} }{dt}\Big\} 
\\&  + \dfrac{\mathrm{I}-\mathrm{sign}\big(A_{i+1/2}\big)}{2}\phi_{i+1/2}+\dfrac{\mathrm{I}+\mathrm{sign}\big(A_{i-1/2}\big)}{2}\phi_{i-1/2} =0
%
\end{split}
\end{equation*}

\begin{figure}[H]
\centering
\subfigure[wave shape]{\includegraphics[width=0.45\textwidth]{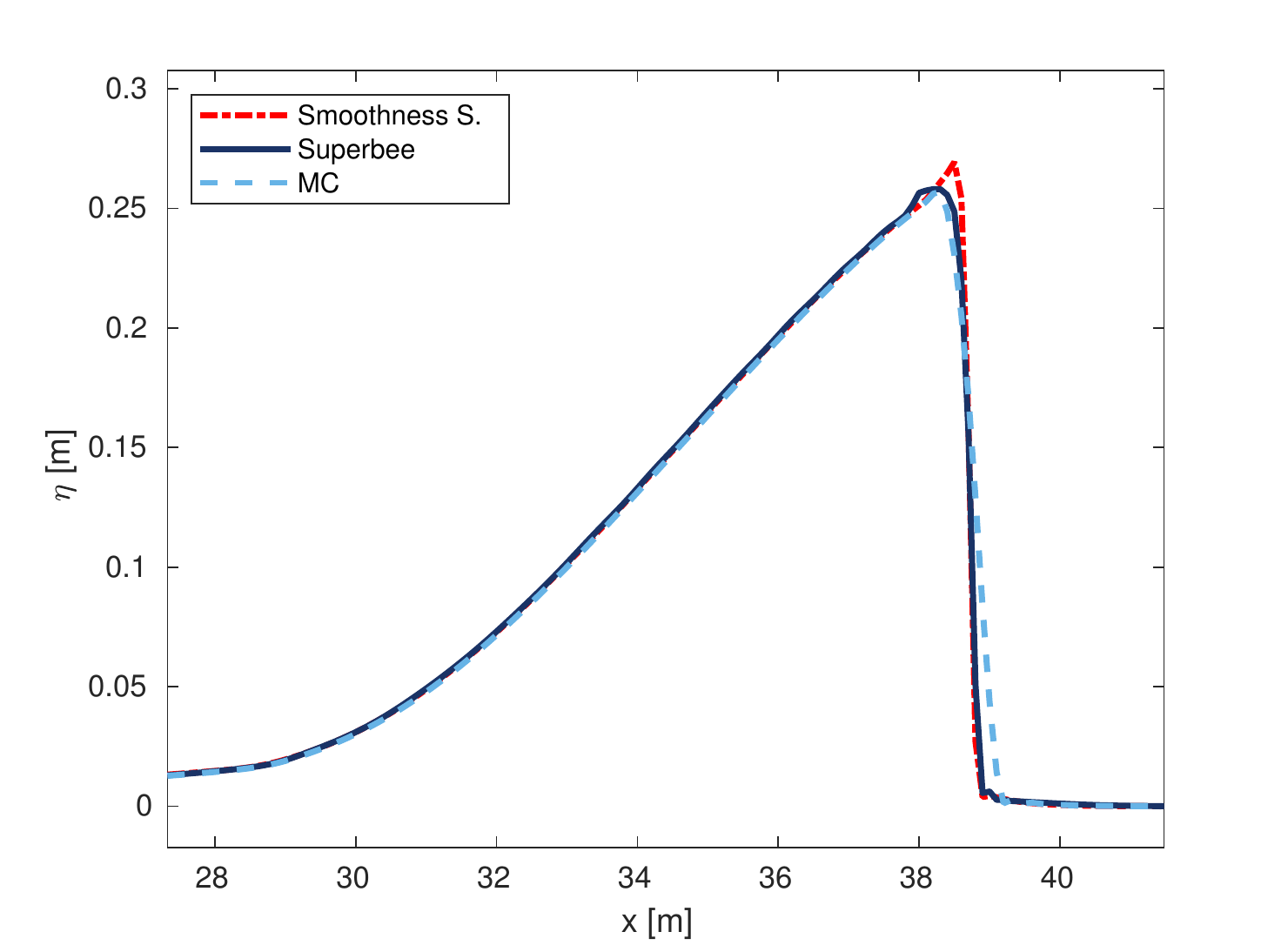}\label{Syno_SW_limit}}
\subfigure[wave shape front zoom]{\includegraphics[width=0.45\textwidth]{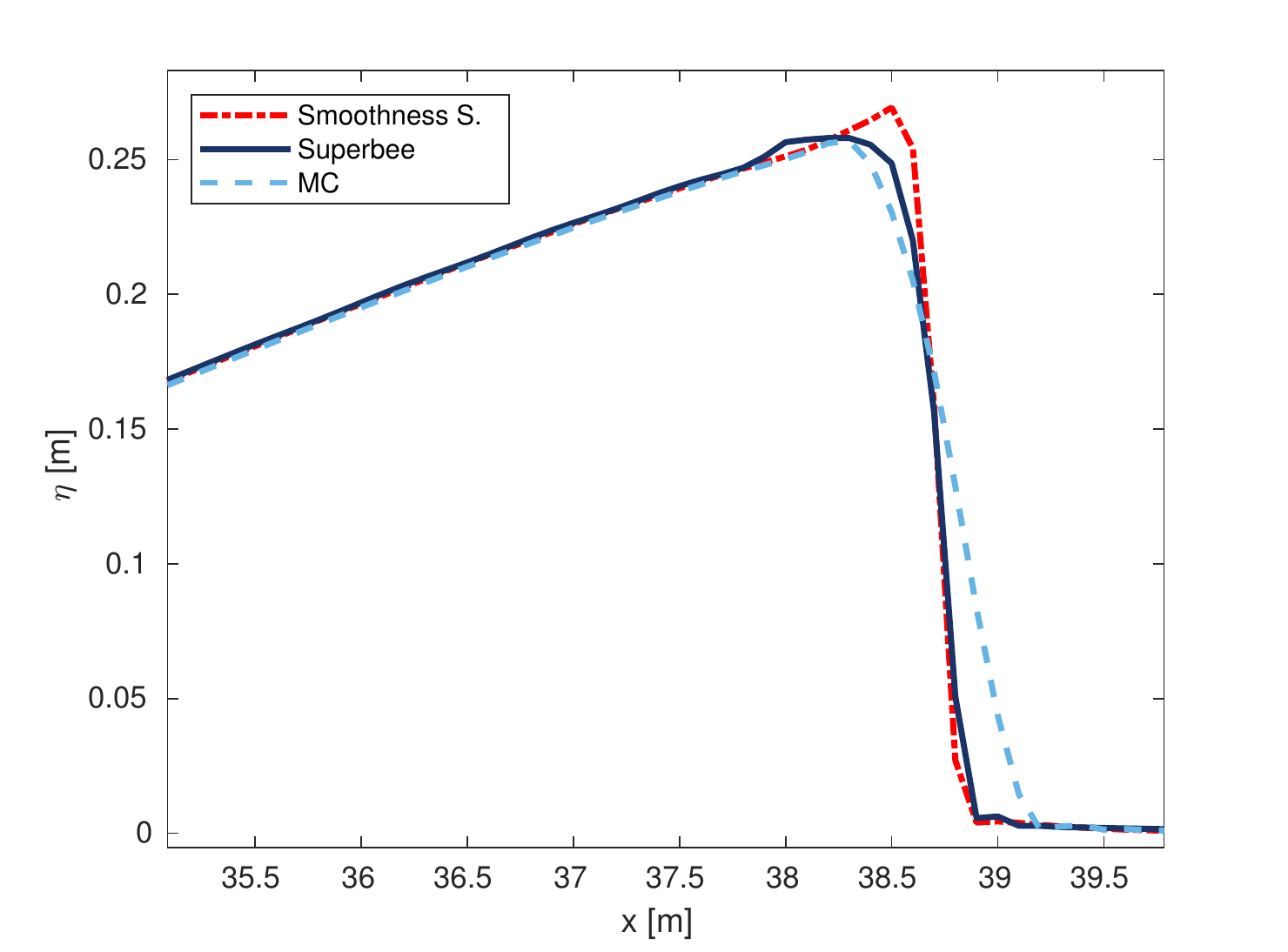}\label{Syno_SW_limit_zoom}}
\caption{Run-up of a solitary wave with $A=0.28$ ffor the pure Shallow Water model with the first-order scheme. Comparison of the wave shapes obtained with different limiters at $t'=15$.}\label{fig_Syno_SW15}
\end{figure}

For the  cell-wise limiter $\delta_{i \pm 1/2}$ we have tested several classical  slope/flux  limiters. In our case, however, it is very important that numerical dissipation, in correspondence of smooth extrema, is not excessive, as this may
hinder the accuracy of the breaking closure. For this reason we have in practice used  a limiter based on a  smoothness indicator, with the main idea to combine two different approximations of a solution derivative. To fix ideas, let us consider the case of a second-order derivative, but for the first-order derivative one proceeds similarly. We consider a fourth-order  approximation
$$
\left(\partial_{xx}\eta\right)^{\text{O4}}=
-\frac{1}{12 \Delta x^2}\eta_{i+2}+ \frac{4}{3 \Delta x^2}\eta_{i+1} - \frac{5}{2 \Delta x^2}\eta_{i} +\frac{4}{3 \Delta x^2}\eta_{i-1} -\frac{1}{12 \Delta x^2}\eta_{i-2}
$$
along with a second-order approximation
$$
\left(\partial_{xx}\eta\right)^{\text{O2}}=  \frac{1}{ \Delta x^2}\eta_{i+1}- \frac{2}{\Delta x^2}\eta_{i}+\frac{1}{\Delta x^2}\eta_{i-1}.
$$
Clearly for a smooth variation of $\eta$ we retrieve
$$
\left(\partial_{xx}\eta\right)^{\text{O4}}= \partial_{xx}\eta +\mathcal{O}(\Delta x^4)\,,\;\;\;\;
\left(\partial_{xx}\eta\right)^{\text{O2}}= \partial_{xx}\eta +\mathcal{O}(\Delta x^2)
$$
and hence 
$$
 \frac{|{\eta_{i+2} - 4\eta_{i+1} + 6\eta_{i} - 4\eta_{i-1} + \eta_{i-2} }|}{12 \Delta x^2}=
 |\left(\partial_{xx}\eta\right)^{\text{O4}}-\left(\partial_{xx}\eta\right)^{\text{O2}}| = \mathcal{O}(\Delta x^2).
$$
On the other hand, for a discontinuity of finite amplitude $\Delta\eta$ around node $i$ we get
$$
\left(\partial_{xx}\eta\right)^{\text{O4}}\approx   \frac{15\Delta \eta}{ 12\Delta x^2}\ \,,\;\;\;\;
\left(\partial_{xx}\eta\right)^{\text{O2}}\approx  \frac{\Delta \eta}{ \Delta x^2},
$$
or, equivalently,
$$
 \frac{|{\eta_{i+2} - 4\eta_{i+1} + 6\eta_{i} - 4\eta_{i-1} + \eta_{i-2} }|}{12 \Delta x^2}=
 |\left(\partial_{xx}\eta\right)^{\text{O4}}-\left(\partial_{xx}\eta\right)^{\text{O2}}| \approx \frac{\Delta \eta}{ 4\Delta x^2}.
$$
We thus define the term
\begin{equation}
\sigma_i=\min(1, \frac{\epsilon + \frac{|\eta_{i}-\eta_{i-1}|}{\Delta x} + \frac{|\eta_{i}-\eta_{i+1}|}{{\Delta x}}}{\epsilon+ \frac{|{\eta_{i+2}- 4\eta_{i+1} + 6\eta_{i} - 4\eta_{i-1} + \eta_{i-2} }|}{12 \Delta x^2} }),
\label{slope1}
\end{equation}
with $\epsilon$ a small parameter used to avoid any division by zero.
In the smooth case, while the numerator in the right slot is bounded, the denominator is of $\mathcal{O}(\Delta x^2)$, and thus the value of $1$ is retained.
Across a discontinuity, the numerator gives $\Delta\eta/\Delta x$, while the denominator gives approximately $\Delta\eta/4\Delta x^2$,
so $\sigma_i = 4\Delta x \ll 1$.  In practice, in every cell we set $\delta_{i\pm1/2} =\min(\delta_i,\delta_{i\pm1})$ where,
 to avoid abrupt variations of the sensor, we set in every node
\begin{equation}
\delta_i=\frac{1}{4}\hat\delta_{i-1}+\frac{1}{2}\hat\delta_{i}+\frac{1}{4}\hat\delta_{i+1}\,,\quad
\hat\delta_i=\left\{
\begin{array}{ll}
\sigma_i\quad&\text{if }\sigma_i < 0.5\\[10pt]
1\quad&\text{otherwise}.
\end{array}
\right.
\label{slope3}
\end{equation}
The advantage of this definition is its ability to better preserve smooth extrema, while still  capturing strong discontinuities.
We have tested this smoothness sensor limiter on several test cases, comparing it to more classical limiters.
We report hereafter the behaviours observed in two typical cases. The first consists in the run-up of a solitary wave, which will be discussed in detail in section \S5.2.
To approximate the solution, we consider the Shallow Water equations with the first-order scheme and test different limiting strategies. 
In Figure \ref{fig_Syno_SW15}, we show the free surface distributions obtained at $t'=15$ (cf. \S5.2 for more details) for the smoothness sensor, the Superbee and the monotonized centered  \cite{LeVeque} limiters. 
Comparing these approximations, one can clearly see that the smoothness limiter improves the preservation of the wave peak.  
\begin{figure}[H]
\centering
\includegraphics[width=0.65\textwidth]{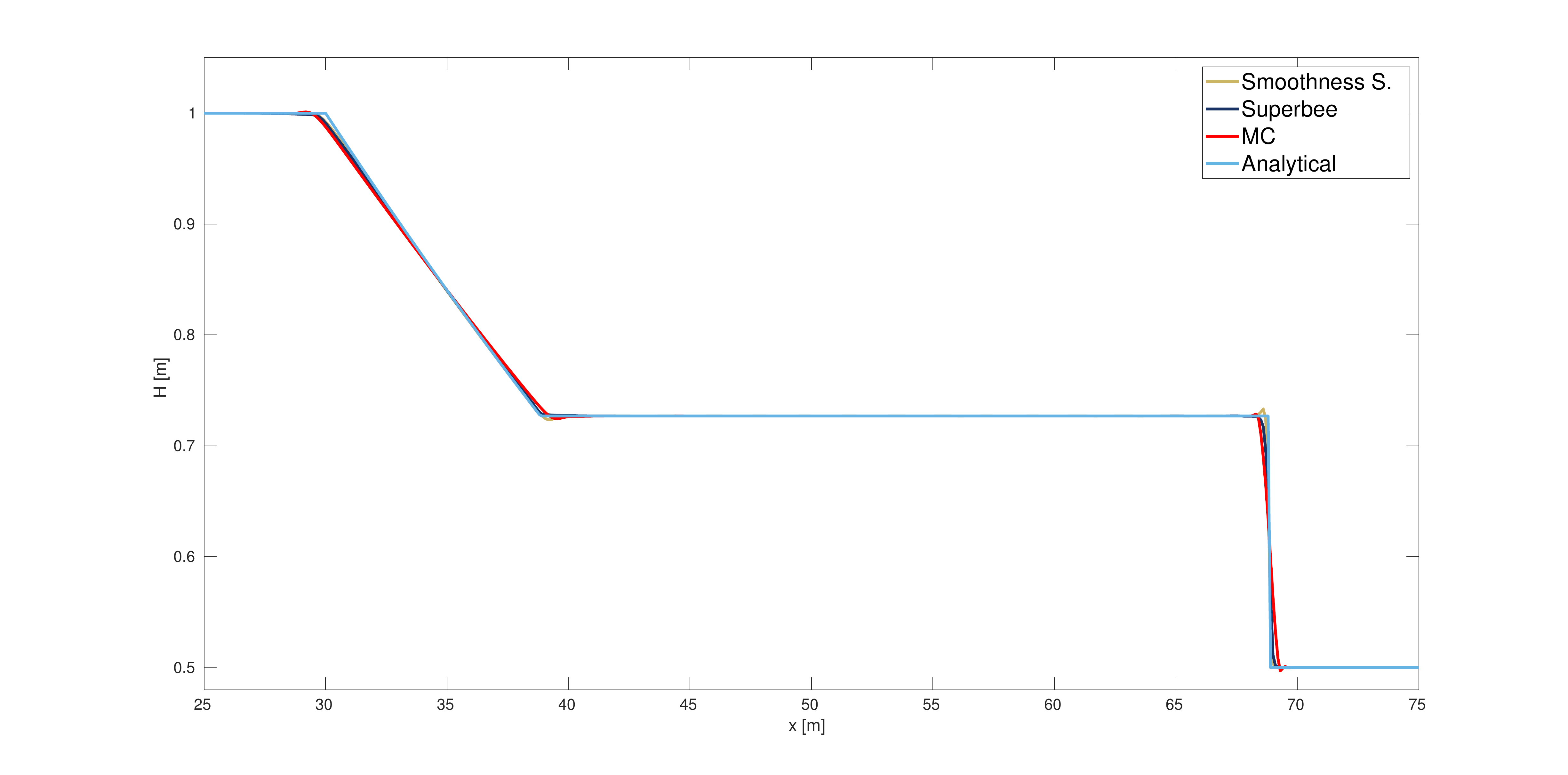}
\caption{Dam break problem. Solution obtained with the smoothness sensor compared to those yield by the Superbee and the monotonized centered (MC) limiters.\label{fig_RP}}
\end{figure}
As a second example, we consider a dam break problem.
The results of Figure \ref{fig_RP} show that the smoothness sensor allows to capture the strong shock, with a behaviour somewhat close to that provided by the monotonized centered limiter.
Note that, to obtain these results, the  upwind flux splitting terms $\text{sign}(A)\Delta F$ need to be   corrected with an entropy fix, for which we have adapted the 
techniques detailed in \cite{Harten,Kermani,Pelanti} (see \cite{bacigaluppi:hal-00990002} for details).
These two benchmark problems allow to highlight the advantage of using the proposed smoothness limiter, and all further test cases will be performed relying on this choice.
In particular,  the complete discretization in its final form reads:
\begin{equation}\label{eq:scheme_full}
\begin{split}
\Delta x & \big(\frac{d\mathbf{W}_i}{dt}  +S_{f_i}\big) + \dfrac{\mathrm{I}-\mathrm{sign}\big(A_{i+1/2}\big)}{2}\phi_{i+1/2}+\dfrac{\mathrm{I}+\mathrm{sign}\big(A_{i-1/2}\big)}{2}\phi_{i-1/2}+ \\
\Delta x&\frac{\mu_{i-\frac{1}{2}}}{2}   \Big\{\frac{1}{3} [\big( \frac{d\mathbf{W}}{dt} +S_f\big)_{i-1} - \big(\frac{d\mathbf{W}}{dt} +S_f\big)_{i}] +  [\text{sign}(A) \big(\dfrac{d\mathbf{W}}{dt}+S_f\big)]_{i-\frac{1}{2}}\Big\}+ \\ 
 \Delta x&\frac{\mu_{i+\frac{1}{2}}}{2}  \Big\{\frac{1}{3} [\big( \frac{d\mathbf{W}}{dt} +S_f\big)_{i+1} - \big(\frac{d\mathbf{W}}{dt} +S_f\big)_{i}]  - [\text{sign}(A) \big(\dfrac{d\mathbf{W}}{dt}+S_f\big)]_{i+\frac{1}{2}}\Big\} =\\
f_{b_i}&\int\limits_{\Omega}\varphi_iD(\mathbf{W}_{\textrm{h}}) dx+ 
 \sum_{i=0}^{N-1} f_{b_{i+\frac{1}{2}}}\int\limits_{x_i}^{x_{i+1}}\left(A(\mathbf{W}_{\textrm{h}})\partial_x\varphi_i\,\tau_{\textrm{s}}\right)  D(\mathbf{W}_{\textrm{h}})dx,
%
\end{split}
\end{equation}
having denoted by $\mu_{i\pm\frac{1}{2}}$ a cell limiter computed as 
$$
\mu_{i\pm\frac{1}{2}} =\left\{ \begin{array}{ll}
1\qquad \quad&\text{if }\;\min(f_{\text{break}_i},\,f_{\text{break}_{i\pm1}}) =1\\
\delta_{i\pm\frac{1}{2}}\omega_{i\pm\frac{1}{2}}\qquad &\text{otherwise}.
\end{array} \right.
$$
Further, $f_{b_i}=f_{\text{break}_i} \omega_i$, where $\omega(H)$ is the exponential wet/dry cut-off introduced in \cite{Ricchiuto2009}, and computed as
$$
\omega_i = \text{exp}\left\{\ -\alpha \Delta x \dfrac{\max\limits_{j=0,N} H_j -\epsilon}{\max( \epsilon, \min\limits_{j=i-2,i+2}H_j-\epsilon )} \right\},
$$
with $\epsilon$ a mesh-dependent small threshold  in practice set to $\epsilon = (\Delta x/L)^2$, where $L$ is  a characteristic length of the domain, and $\alpha$ is 
usually set to $\alpha\approx 10$. In practice we set  $f_{b_{i\pm1/2}}=\min(f_{b_{i\pm 1}},f_{b_{i}} )$ and similarly for $\omega$.\\
The resulting scheme can be readily shown to be well-balanced in wet regions. In particular considering the lake at rest state given by  $q=0$ and $\eta$ constant across the entire domain,
one can show that all the dispersive terms, as well as the hydrostatic fluctuations $\phi_{i\pm1/2}$ and the friction terms,  are identically zero.  
 The scheme thus reduces to $d\mathbf{W}/dt =0$, which shows that this state is indefinitely preserved.
More details of the analysis can be found in \cite{rf14}.
In this work we have generalized the validity of this property, to account also for dry areas, adopting the modifications  proposed in \cite{Castro2005,Delis2008,Ricchiuto2009}.
Finally, the integration in time of  \eqref{eq:scheme_full} is performed with a non-dissipative implicit Crank-Nicholson method, and 
wave generation and boundary conditions are handled using the  source functions and absorbing layers discussed to some extent in \cite{rf14}.
The interested reader can refer to this reference, and to  \cite{bacigaluppi:hal-00990002}, and references therein, for more details.

\paragraph{\textit{Remark 1}}
\revP{Note that the choice of the presented method is not crucial. However, its detailed description allows to underline the strict interplay between the numerics (limiting procedure), and the modelling of wave breaking. While other schemes have been used in combination with a similar wave breaking closure (cf. \cite{Kazolea14,filippini2016,kazolea2018}),  
such interaction is always necessary. This  makes the understanding of the stabilisation and shock-capturing mechanisms critical.}
\paragraph{\textit{Remark 2}}
\revP{While the numerical approximation has been extensively tested for convergence on the Madsen and S{\o}resen equations in \cite{rf14}, our focus here shall be on the breaking and shock capturing features.}
\section{Wave Breaking Detection and Implementation}
\label{sec:wavebreak}

We discuss in this section the detection criteria considered in this paper.  The objective is to provide computable definitions for the breaking flag $f_{\text{break}}$ required in the scheme.
In particular, we consider three  approaches: one based on wave non-linearities, taken from \cite{Tonelli2011,Tonelli2012}; the second defined on both a slope and a vertical velocity condition, taken from
\cite{Kazolea14}; the third being our implementation of the so-called convective or critical free surface Froude criteria discussed e.g. in  \cite{LH69} (see also \cite{Bjorkavag2011,Brun2018} and references therein). 
To our knowledge this is a first attempt at using such a convective criteria in predicting wave breaking in the context of Boussinesq simulations with a breaking closure.

From a physical point of view,  it is well known that one of the main parameters ruling wave breaking is the Froude number.
However, when looking at wave breaking, several  definitions  of this quantity are possible, going from the classical local value $Fr =|u|/\sqrt{gH}$, to definitions based on the ratio of the wave amplitude over the depth (cf. \cite{Tonelli2011}, or \cite{chanson2004} for open channel flows), to
other definitions involving wave by wave analysis as the local trough Froude criteria proposed in \cite{Okamoto2006}, and many others.

For moving bores,  the well known studies by Favre \cite{favre1935}  and Treske \cite{treske1994} show that 
the limiting value of the local Froude number, relative to the bore above which breaking occurs, is about $1.3$.  
In general, the onset of breaking  has been identified to happen in a range of the Froude numbers between $1.3$ and $1.6$ (see \cite{Okamoto2006,Tonelli2009,Tonelli2010,Tonelli2011,Tonelli2012}). 
In practice, defining a local Froude number starting from  data issued from numerical (or experimental) data is not a trivial task as many parameters are less clearly identified than, e.g.,
in the experimental set-up of  \cite{favre1935}  and  \cite{treske1994}. In a wave by wave analysis, for example, there is the question of determining the wave's celerity, the location in which the velocity and Froude condition
 should be estimated, which velocity should be used, and so on. 
 
The three criteria considered here allow to cover good  part of the existing detection methods.
To be able to fairly compare them,  our objective is to provide approximations to the sought solution, based on the same numerical method, and within the same breaking closure approach.  
As much as possible of their practical implementation is also discussed, 	to enable for reproducibility.


\subsection{A Criterion Based on Local Non-linearity}
\label{sec:Tonelli} 
The  first criterion we have considered is the simplest to set up. It  is taken from \cite{Tonelli2009,Tonelli2010,Tonelli2011,Tonelli2012} in which the authors 
 present a detection method based on  the non-linearity parameter $\mathcal{E}$, defined as the ratio wave amplitude over still water depth, i.e.
\begin{equation}
\mathcal{E}=\frac{A}{h}.
\label{epsilon1}
\end{equation}
For the onset of  breaking,  the value $\mathcal{E}= 0.8$ is suggested. This is  obtained assuming a transition at  a Froude $\approx 1.6$ \cite{Okamoto2006,Tonelli2009,Tonelli2012},
 and considering that  the  limit value  to obtain a  stable solitary wave is $\mathcal{E}\approx 0.78$ \cite{Kazolea2013}.
Clearly, this criterion is purely local, as  it involves no variations (derivatives) of the wave characteristics. 
Concerning the termination criterion, once $\mathcal{E}>0.8$, the breaking closure is kept active until $\mathcal{E}$  is back within a range between $0.55 \div 0.25$, as suggested in  \cite{Tonelli2011,Tonelli2012}. \\

We will refer in the following to this detection method as the  \textbf{Local Criterion}. It has the advantage of being simple, but in some cases it 
lacks generality, as, for example,  in presence of wave set-up (or set-down), and in general 
when the  value of $h$ is not anymore relevant  to define a reference depth or wave celerity.\\

In our implementation of this detection method,  we have computed in   each node   the discrete analogue of   \eqref{epsilon1} given by
\begin{equation}
\mathcal{E}_{\textrm{h}}=\frac{|\eta|}{|h|+\epsilon}\quad ,
\label{epsilonimplem}
\end{equation}
with $\epsilon\approx (\Delta x/L)^2$, which allows to avoid any division by zero when considering the run-up on sloping beaches.
Further, whenever a point in space satisfies $\mathcal{E} > \mathcal{E}_{\text{critical}}$, with $\mathcal{E}_{\text{critical}}=0.8$,
we seek the point downstream in the wave direction satisfying $\mathcal{E}=0.3$. The downstream direction is  either known a priori,
or determined using the technique discussed  in section \S4.3.  In the whole region between the two points we set $f_{\text{break}}=0$.

\subsection{Hybrid Slope/Vertical Velocity Criterion}
\label{sec:Delis} 
The second criterion we consider was initially proposed in  \cite{Kazolea14}. It combines two criteria: one related to  the steepness of the wave \cite{Schaffer1993} and
the second involving the vertical displacement of the free surface \cite{bfwts06}.
 In practice, the onset of breaking is determined by the  validity of one of the two following conditions
\begin{equation}
\begin{split}
& |\partial_x \eta| \geq \tan \phi \quad  \text{with}\quad \phi\in[14^{\circ},33^{\circ}]\\[5pt]
& \dfrac{\partial_t \eta}{\sqrt{g |h|}} >\gamma\quad \text{with} \quad 0.35\leq \gamma \leq 0.65.
\end{split}
\label{secondcritI}
\end{equation}
In \cite{Kazolea14}, the value  $\gamma=0.6$ is prescribed as a default.  Note that the second condition closely resembles a Froude number condition based on the velocity of the free surface.
This condition is inefficient for steady hydraulic jumps, in which case the slope criterion is more effective.
For this the choice of the critical angle  depends strongly on the breaking type. As in \cite{Kazolea14} and \cite{Tissier2012}, and unless explicitly mentioned,
 we have used the default value  $\phi=30^{\circ}$.  The termination of breaking is based in this case on the Froude condition proposed in \cite{Tissier2012}, which relies on the analysis of moving bores.
  Following \cite{Kazolea14}, this detection method is referred to as  \textbf{Hybrid Criterion}.

\begin{figure}[H]
\centering
\includegraphics[scale=0.45]{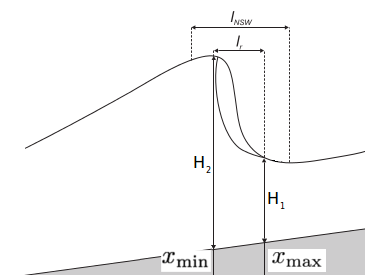}
\caption{Sketch of the lengths that determine the roller magnitude.\label{fig_lNSW}}
\end{figure}

The practical implementation   is done in four phases as follows.
\begin{enumerate}
\item \textit{Pre-flagging.} We compute in every node
\begin{equation*}
\begin{split}
&Fr_{w_i}^n:=\dfrac{|\eta_{i}^{n} -\eta_i^{n-1}|}{\Delta t\sqrt{g( |h_i| + \epsilon)}}\,\quad \text{and}\\ &  \nabla\eta_i^n :=\max\left( \dfrac{|\eta_{i+1}^n-\eta_{i-1}^n|}{2\Delta x},\,\dfrac{|\eta_i^n-\eta_{i-1}^n|}{\Delta x},\,\dfrac{|\eta_{i+1}^n-\eta_i^n|}{\Delta x} \right).
\end{split}
\end{equation*}
If 
$$
Fr_{w_i}^n > \gamma \quad\text{or}\quad \nabla\eta_i^n \ge \tan\phi,
$$
we set $f_{\text{break}_i}=0$. Adjacent cells involving flagged nodes are clustered, thus obtaining a preliminary definition of the breaking length.
\item  \textit{Peak and trough water level.} In each cluster of breaking cells we determine the maximum and minimum values of the free surface elevation,  $H_2$ and $H_1$, and the corresponding mesh points $x_{\max}$ and  $x_{\min}$, as  on  Figure \ref{fig_lNSW}.
\item \textit{Termination check.} We use the criteria proposed in  \cite{Tissier2012} and compute the bore Froude number  $Fr_b$ defined as 
\begin{equation}
Fr_b=\sqrt{\frac{\big(2\frac{H_2}{H_1}+1\big)^2-1}{8}}.
\label{froudeII}
\end{equation}
If $Fr_b <Fr_{b_{\text{cr}}}$, we set  $f_{\text{break}}=1$  along all the breaking length and move to the next cluster of breaking cells. We set   $Fr_{b_{\text{cr}}}=1.3$, as proposed  in  \cite{Kazolea14,Tissier2012}. 

\item \textit{Roller magnitude estimation and final flagging.} If $Fr_b >Fr_{b_{\text{cr}}}$,  we follow again the definitions of \cite{Kazolea14,Tissier2012}: 
we compute the physical roller length as $l_r=2.9(H_2-H_1)$, and we define a  numerical roller of length $l_{NSW}=2.5 l_r$.
We set $f_{\text{break}_i}=0$  in all points verifying $ 2 x_i -(x_{\max}+x_{\min}) < \bar{l}_{NSW} $.  
\end{enumerate}

\subsection{The Convective or Critical Free Surface Froude Criterion}
\label{sec:Philippe} 
The main idea behind this detection method is that wave breaking occurs when the free surface velocity at the  peak of the wave exceeds the celerity of the wave.
Therefore,  one may  consider the  Froude number
\begin{equation}
Fr_s=\frac{u_s}{c_b},
\label{FroudeCrit}
\end{equation}
with $u_s$ the free surface velocity at the wave crest, and $c_b$ the  celerity.
If $Fr>1$ then the wave should break.  In the following we will refer to this detection method  as to the 
\textbf{Physical Criterion}.
 
To implement this method in practice, we must provide computable definitions of  both the free surface velocity and the celerity.
Concerning the first, we recall the  asymptotic development of the horizontal velocity reading  \cite{book_lannes}
$$
u(z) = u  -  \Big[\dfrac{(z-h)^2}{2} -   \dfrac{H^2}{6} \Big] \partial_{xx}u 
-  \Big( \dfrac{H}{2}+z-h \Big)(2\partial_xh\partial_x u+\partial_{xx}hu),
$$
where $u$ is the depth averaged velocity. The idea is then to evaluate the above expression on the free surface, for $z=\eta=H+h$,
and at the wave's crest, where it is assumed that $\partial_xu=0$. Evaluating the above expression,  we obtain
\begin{equation}\label{eq:us}
u_s = u  -   \dfrac{H^2}{3} \partial_{xx}u  - \dfrac{3q}{2} \partial_{xx}h.
\end{equation}
Concerning the  celerity of the wave,  we combine the two relations
$$
\partial_t\eta+c_b\partial_x \eta =0 \quad\text{and}\quad\partial_t\eta+ \partial_x q =0,
$$
obtaining 
\begin{equation}\label{eq:cb}
c_b =   \dfrac{\partial_xq}{\partial_x\eta} \,.
\end{equation}
With these definitions, we proceed to an implementation based on the following four steps.
\begin{enumerate}
\item \textit{Pre-flagging.} Same as for the hybrid criterion (see section \S4.2).
\item  \textit{Peak and trough water level.} Same as for the hybrid criterion (section \S4.2).
\item \textit{Evaluation of $Fr_s$.} We need to evaluate the free surface velocity in correspondence of $x_{\max}$, as well as the celerity.
For the discrete free surface velocity $u_{s_{\textrm{h}}}$,  we use  \eqref{eq:us} in which  the second derivative of the velocity is computed first by means of 
 standard second-order centred finite difference formula, and then smoothed, in order to remove some of the numerical noise. The filtering is done using one iteration of a 
5-point moving average, using 2 mesh points on either side of $x_{\max}$. We neglect the second term in  \eqref{eq:us} in all the computations performed, which involve piecewise constant slopes. The implementation of \eqref{eq:cb}, used to evaluate the celerity, is 
\begin{equation}
c_{b_\textrm{h}}= \dfrac{q(x_{\max}) -q(x_{\min})  }{\eta(x_{\max}) - \eta(x_{\min}) + \epsilon},
\label{cb1}
\end{equation}
with $\epsilon$ added to avoid any division by zero. The critical free surface Froude is then discretely evaluated as 
\begin{equation}\label{eq:Frh}
Fr_{s_{\textrm{h}}}=\dfrac{|u_{s_{\textrm{h}}}|}{|c_{b_\textrm{h}}|+\epsilon}\,.
\end{equation}

\item \textit{Roller magnitude estimation and final flagging.} If $Fr_{s_{\textrm{h}}} >Fr_{s_{\text{cr}}}$,  we  proceed as in section \S4.2:
we define  $l_{NSW}=2.5 l_r$ with $l_r=2.9(H_2-H_1)$, and 
we set $f_{\text{break}_i}=0$  in all points verifying $ 2 x_i -(x_{\max}+x_{\min}) < \bar{l}_{NSW} $.  
The critical value by default used is $Fr_{s_{\text{cr}}}=1$, but, as we shall see in the numerical experiments, other values have been tested to compensate for the shoaling behaviour of the considered Boussinesq model.
\end{enumerate}

\noindent \paragraph{\textit{Remark}} \revP{The determination of the wave celerities necessary for the convective criterion is a key point for future
extensions. This issue is common to both numerical and experimental techniques aiming at characterising in detail the wave-field. The filtering approach used here is a first solution, but more involved techniques will be for sure required to implement this approach in a more general context, both concerning the quality of the estimated celerities and the extension to multiple space dimensions.}

\subsection{Additional Implementation Details}

As discussed in detail in \cite{kazolea2018}, the wave breaking closure used here introduces easily oscillations at the interface of the coupling between the Bousisnesq model and the Shallow Water equations.
Indeed, in \cite{kazolea2018}, it is shown that this is especially important when a mesh refinement is involved, as, if uncontrolled, it would lead to the blow up of the computation, instead of the sought mesh convergence.  In particular, on coarse meshes the amplitude of these perturbations remains very small.
However, since there is no clear understanding of the link between the mesh size and numerical spurious effects, we have tried to minimize the sources of instability. 
To this end, some techniques for the minimization of the oscillations have been already mentioned, as, e.g., the use of filtering to evaluate the second derivatives of the velocity.  Furthermore, through the determination of the roller magnitude, as, e.g. within the fourth step of the practical implementation of the physical criterion, we have also avoided the introduction of spurious breaking. Specifically, we have throughout avoided the  introduction of rollers with a size smaller than the stencil of our numerical scheme, i.e. for $\bar{l}_{NSW}\le 4\Delta x$,  we set back $f_{\text{break}_i}=1$.
Similarly, if two breaking regions are too close, they are merged, such that, if  the distance between two breakers is less than   $ 4\Delta x$, we set  $f_{\text{break}_i}=0$  also for the nodes in between.

Finally, we avoid any possible wrong detection, as, e.g., the  breaking  being detected  on the rear side of the wave via the check of the fulfilment of $c_{b_\textrm{h}}(x_{\max}-x_{\min})\ge 0$.

\section{Numerical Tests and Results}
\label{cha:5numerical}
\subsection{Wei's Solitary Wave Shoaling on a Slope}
\label{WeiTest}
In order to study the physics of the wave breaking,
we consider the benchmark problems proposed in Wei et al. (\cite{Wei95}), 
involving  solitary waves shoaling over different constant slopes. 
 These tests are used to investigate the physical behaviour obtained with the different breaking criteria.

\subsubsection{Case 1:  Slope of 1:35}
\label{WeiTest35}
The set-up is given by a solitary wave of amplitude $A = 0.2$m over a spatial domain defined within $[-12, 40]$m. The  bathymetry, set initially with $h_0 = 1$m
depth, has a constant slope 1:35 between $0 \leq x \leq 35$m and continues till the end of the domain with a flat slope. The solitary wave starts at the toe of the slope set at $x = 0$m. The total considered time is $30$s.
The time step is dictated by the CFL number $\nu= 0.3$ on a mesh with $\Delta x=0.05$m. 
Absorbing boundary conditions are set on the left and right ends of the domain using sponge layers $3$m wide (see \cite{rf14} for details).
The wave's evolution is described as a function of a dimensionless time $t' = t \sqrt{g/h_{0}}$.

\begin{figure}[h]
\centering
\subfigure[Local]{\includegraphics[width=0.71\textwidth]{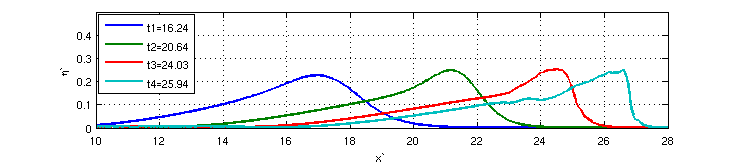}\label{Wei35_local}}
\subfigure[Hybrid ]{\includegraphics[width=0.71\textwidth]{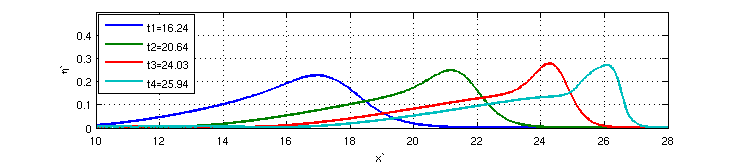}\label{Wei35_Hybrid}}
\subfigure[Physical  ($Fr_{s_{\text{cr}}}=1$)]{\includegraphics[width=0.71\textwidth]{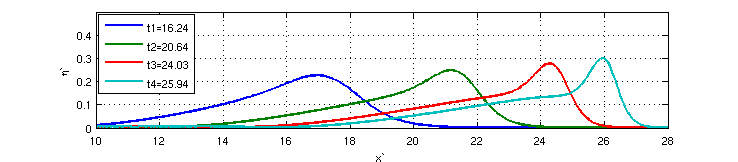}\label{Wei35_off}}
\caption{Benchmark problem of \cite{Wei95} with slope 1:35. Comparison between different wave breaking criteria.}\label{fig_Wei35}
\end{figure}

As in   \cite{Wei95}, we will  focus on  four different discrete  times: $t_1' = 16.24$, $t_2' = 20.64$, $t_3' = 24.03$ and $t_4' = 25.94$. 
We visualize the wave profile at these times when using the three different detection criteria: local (Figure \ref{Wei35_local}), hybrid (Figure  \ref{Wei35_Hybrid}), and  physical
(Figures \ref{Wei35_off}).   
\vspace{-0.3cm}
\begin{figure}[h]
\centering
\subfigure[Local ]{\includegraphics[width=0.47\textwidth]{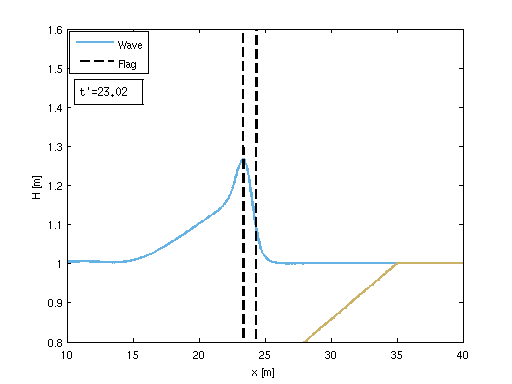}\label{35B1}}
\subfigure[Hybrid ]{\includegraphics[width=0.47\textwidth]{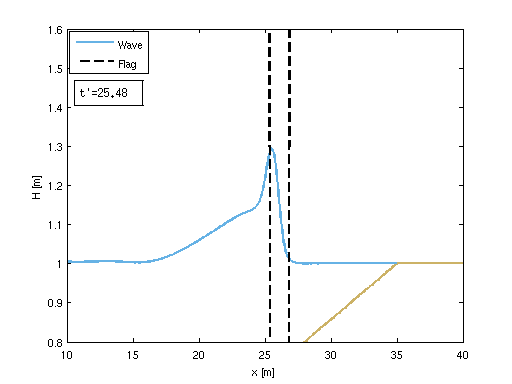}\label{35B2}}
\subfigure[Physical ($Fr_{s_{\text{cr}}}=1$)]{\includegraphics[width=0.47\textwidth]{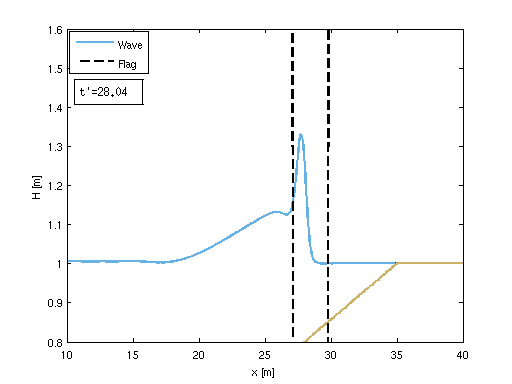}\label{35B4}}
\subfigure[Physical  ($Fr_{s_{\text{cr}}}=0.75$)]{\includegraphics[width=0.47\textwidth]{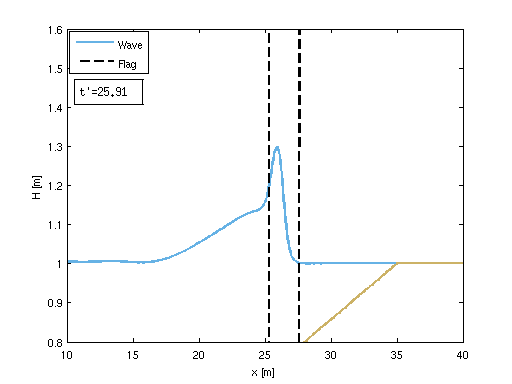}\label{35B3}}
\caption{Benchmark problem of \cite{Wei95} with slope 1:35. First instants of wave breaking activation for different wave breaking detection criteria.}\label{fig_Wei35break}
\end{figure}
From the figures we can see that all criteria provide the same wave profile at $t_1'$ and $t_2'$, at which the wave propagates and starts shoaling, becoming taller and less symmetric.
The shoaling continues at   $t_3'$ with the hybrid and physical criteria, but 
with the local criterion the smooth wave peak has been already flattened, which is the sign of an early onset of wave breaking.
This becomes even clearer  at time $t_4'$, at which this criterion  provides a wave shape very close to that of a shock. 
 The hybrid criterion fares better,  however, the first signs of an early activation of the breaking closure are outlined.
 Indeed, the wave at $t_3'$  is still smooth, however we can see the wave height reducing from  $t_4'$ can be seen as the wave profile is already shorter
 than the previous one, a sign that breaking has taken over the shoaling phase. 
The physical criterion, instead, shows wave profiles still growing from  $t_3'$  to $t_4'$, meaning that  shoaling  is  still the dominating phenomenon. 
As a general remark,  compared to the plots in  \cite{Wei95}, the   wave shapes are  shorter and less steep.
This is related to the fact that  the  Boussinesq model \eqref{eq:Bsqsys_eq}  under-shoals compared both to fully non-linear models and to
other weakly non-linear models known to be prone to over-shoaling, as, e.g.,  the one used in the reference. 
This is a known phenomenon. The interested reader can refer to  \cite{fbcr15},
 where the behaviour of well-known weakly non-linear Boussinesq models, including the one used here,
is studied in the transition between weakly and fully non-linear waves.
To get further insight into the behaviour of the different solutions,   we have singled out the first time of breaking for all the considered detection criteria.  
In Figure  \ref{fig_Wei35break}, we visualize the wave's position and shape, as well as the closure region at the first breaking instant computed with the procedure described in section \S4.
As shown in Figure \ref{35B1}, the local criteria provides a very early first onset of wave breaking at $t'=23.02$. The figure also shows that the breaking termination condition used within this criteria gives a very short sized roller, which does not reach the foot of the wave.
Figure \ref{35B2}  provides the same visualization for the hybrid criteria. The onset of breaking is at  $t'=25.48$, much closer to the reference one. The  size of the roller is, in this case, larger, and allows to cover the whole front of the wave, including its foot. This is most likely a better situation for the stability of the criteria, as the steeper parts of the wave  are all embedded in the roller. Passing to the physical criterion,  we can see in Figure  \ref{35B4} that for $Fr_{s_{\text{cr}}}=1$, this method provides a late onset of breaking at $t'=28.04$.
This is a consequence of the undershoaling properties of the Madsen and S{\o}rensen model, that provide waves with reduced steepness, and thus curvature at the peak. This leads to a late explosion of the free surface velocity $u_s$ in \eqref{eq:us}. The lateness in the onset of the breaking also leads to 
a considerably taller wave, and, as a consequence, a roller thickness much larger than that observed, e.g., in Figure  \ref{35B2}. To compensate for the behaviour of the Boussinesq model, several solutions are possible. One could be, for example, to modify \eqref{eq:us} in order to compensate for the aforementioned properties of the model, e.g. by replacing the vertical development  $u(z)$ by a development on $q(z)$, using ideas discussed in  \cite{fbcr15}. Here we adopt a simpler  solution, and introduce a reduced value of the critical Froude, set to $Fr_{s_{\text{cr}}}=0.75$. 
\vspace{-0.4cm}
\begin{figure}[h]
\centering
\subfigure[Local ]{\includegraphics[width=0.71\textwidth]{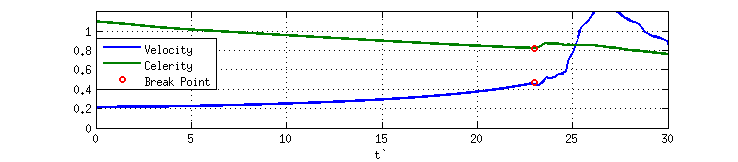}\label{UC35_L}}
\subfigure[Hybrid ]{\includegraphics[width=0.71\textwidth]{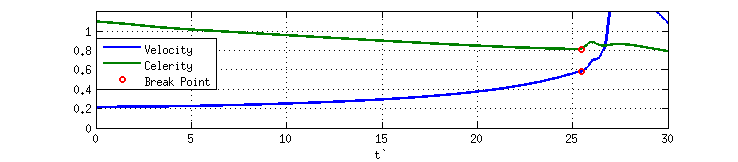}\label{UC35_H}}
\subfigure[Physical  ($Fr_{s_{\text{cr}}}=1$)]{\includegraphics[width=0.71\textwidth]{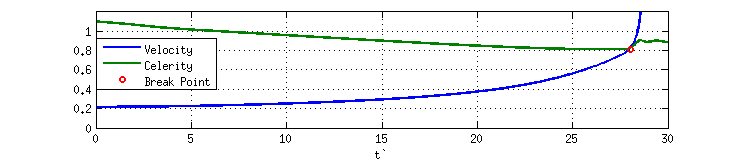}\label{UC35_1}}
\subfigure[Physical  ($Fr_{s_{\text{cr}}}=0.75$)]{\includegraphics[width=0.71\textwidth]{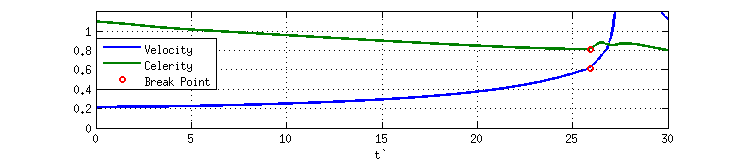}\label{UC35_75}}
\caption{Benchmark problem of \cite{Wei95} with slope 1:35. Parametric  study of the breaking criteria: time evolution of the celerity and of the  free surface  velocity at the peak. }\label{fig_UC35}
\end{figure} 
As shown in Figure \ref{35B3} this provides a breaking time $t'=25.91$,   within $0.1\%$ of the reference, and  a slightly reduced roller thickness.
The wave profiles obtained in this case are identical to those of Figure \ref{Wei35_off}.  
In  \cite{Wei95} a very interesting analysis is reported visualizing the dynamics of the wave 
in terms of evolution of the peaks celerity versus the particle speed at the peak. This is particularly relevant here,
as the physical criterion is actually based on the location where these two curves cross each other.
 In Figure \ref{fig_UC35} we compare  the temporal evolution of the non-dimensional wave's celerity $c_b'= c_b/\sqrt{gh_0}$ to the non-dimensional particle velocity of the wave's crest $u_s'=u_s/\sqrt{gh_0}$.
Note that in the reference the authors report the celerity obtained from  $ \hat{c}_w' =\frac{dx_c}{dt}/\sqrt{gh_0}$, where $x_c$ corresponds to the position of the wave's crest.
Here, we use both for  $c_b$ and for $u_s$ the discrete definitions proposed    in section \S4.3, so that the breaking onset for the physical criterion is exactly represented on the plot.\\
 In figure \ref{fig_UC35}, we also represent with   red circles  the first activation of breaking (represented as a circle in 
 the reference  as well).
 Even though neither of the underlying models, nor the computation of the celerity and of the particle speed are the same as in  \cite{Wei95},
the behaviours observed are  very close to those of  the reference.\\
As the wave shoals, with all the criteria  we obtain an increase in the particle velocity at the crest,
 and a decreasing  celerity at the front. This is  expected for shoaling conditions. The results allow to see graphically the early onset of breaking for 
 some of the criteria, which is, of course, natural as the only one actually using the physical quantities involved is the physical criterion.
It is interesting to note, however, that after the start of the  breaking, all closure models provide a  strong increase in the particle velocity at the crest, which is absent in the reference. 
This is certainly a consequence of the wave profile sharpening, and eventually converging into a shock, 
and of the use of \eqref{eq:us}  to  compute $u_s$.  Interestingly, very similar excursion in the free surface velocities in breaking waves can be  
observed experimentally, see e.g. \cite{Melville1988}.

\subsubsection{Case 2:  Slope of 1:15}
\label{WeiTest15}
\begin{figure}[H]
\centering
\subfigure[Physical  ($Fr_{s_{\text{cr}}}=1$)]{\hspace{0.43cm}\includegraphics[width=0.71\textwidth]{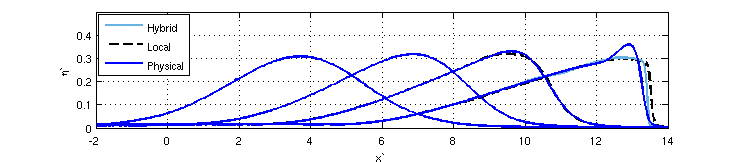} \label{Complete1a}}
\subfigure[Physical  ($Fr_{s_{\text{cr}}}=0.75$)]{\hspace{0.43cm}\includegraphics[width=0.71\textwidth]{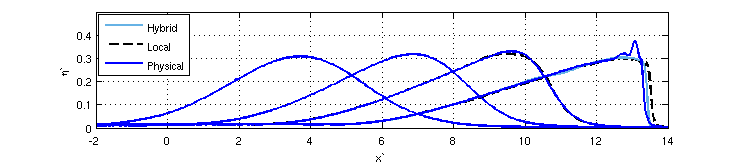} \label{Completea}}
\caption{Benchmark problem of \cite{Wei95} with slope 1:15. Direct comparison of the wave profiles for two different choices of $Fr_{s_{\text{cr}}}$ within the Physical Criterion.}\label{fig_WeiComplete1b}
\end{figure}

\begin{figure}[H]
\centering
\subfigure[Local ]{\includegraphics[width=0.715\textwidth]{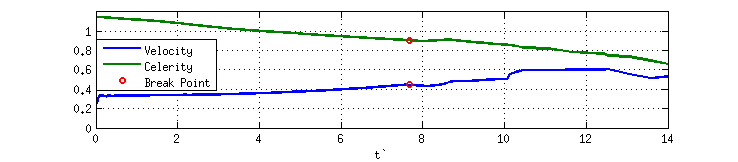}\label{UC15_L}}
\subfigure[Hybrid ]{\includegraphics[width=0.715\textwidth]{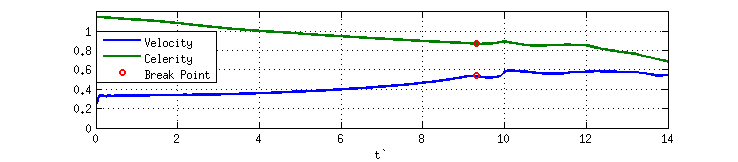}\label{UC15_H}}
\subfigure[Physical  ($Fr_{s_{\text{cr}}}=1$)]{\includegraphics[width=0.715\textwidth]{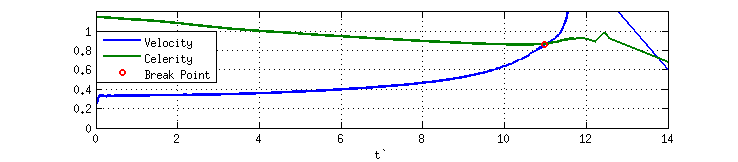}\label{UC15_1}}
\subfigure[Physical  ($Fr_{s_{\text{cr}}}=0.75$)]{\includegraphics[width=0.715\textwidth]{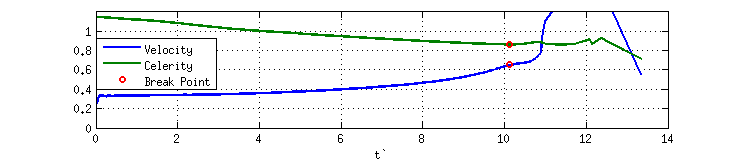}\label{UC15_off}}
\caption{Benchmark problem of \cite{Wei95} with slope 1:15. Parametrical study of the breaking criteria: time evolution of the celerity and of  the   free surface  velocity at the peak.}\label{fig_UC15Ph}
\end{figure}

We consider the same problem set-up,  with a solitary wave of amplitude  $A=0.3$m, shoaling on a slope 1:15.
 Following \cite{Wei95}, we  start looking at the dimensionless time instants $t_1'=3.23$, $t_2'=6$, $t_3'=8.4$ and $t_4'=11.32$. Differently from the previous case,
 these instants do not include a breaking point for the reference \cite{Wei95}.   In our case, all criteria provide, nevertheless, breaking before $t_4'$, as shown in Figure \ref{fig_WeiComplete1b}. 
\vspace{-0.4cm}
    \begin{figure}[H]
\centering
\subfigure[Local ]{\includegraphics[width=0.47\textwidth]{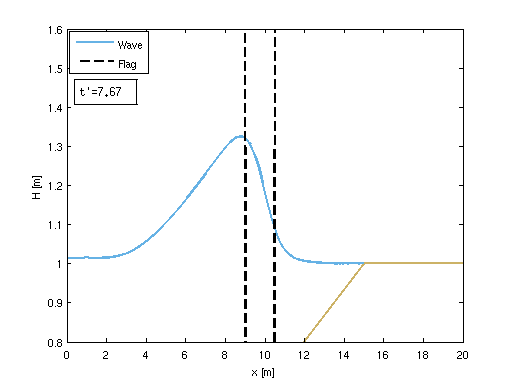}\label{15B1}}
\subfigure[Hybrid ]{\includegraphics[width=0.47\textwidth]{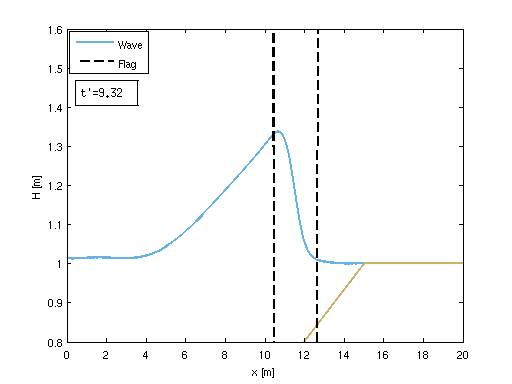}\label{15B2}}
\subfigure[Physical  ($Fr_{s_{\text{cr}}}=1$)]{\includegraphics[width=0.47\textwidth]{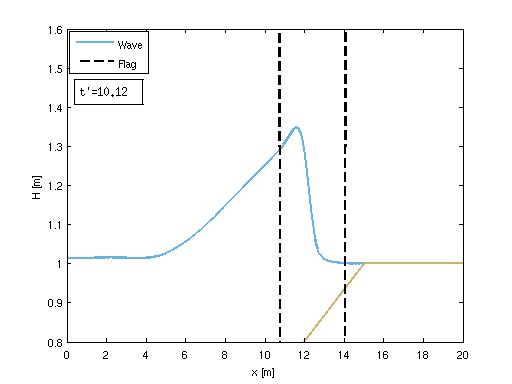}\label{15B3}}
\subfigure[Physical  ($Fr_{s_{\text{cr}}}=0.75$)]{\includegraphics[width=0.47\textwidth]{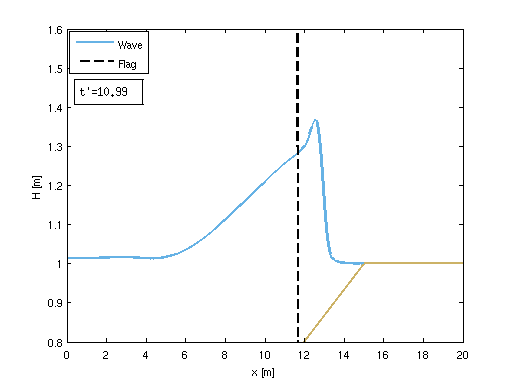}\label{15B4}}
\caption{Benchmark problem of \cite{Wei95} with slope 1:15 - Comparison between different wave breaking criteria.}\label{fig_Wei15break}
\end{figure}
 As before, the earliest onset is obtained with  
  the local criterion for  $t'=7.67$.
 The hybrid detection method predicts a first breaking at $t'=9.32$, while the physical 
 criterion gives breaking onset times of  $t'=10.99$ and  $t'=10.12$ with
    $Fr_{s_{\text{cr}}}=1$  and   $Fr_{s_{\text{cr}}}=0.75$ respectively. The latter  thus provides the latest breaking onset, and gives 
    wave profiles  very close to the ones reported in \cite{Wei95},     including the peak within the last wave which has not dissipated yet.
       These observations are confirmed by the wave profiles at the first breaking instant, reported in  Figure \ref{fig_Wei15break}.
    As in the previous case, we see that the local criterion gives a smaller roller which does not involve the foot of the wave, while
    the roller  obtained with the hybrid criterion encompasses the whole front of the wave. The physical criterion with $Fr_{s_{\text{cr}}}=1$
    provides a taller wave with a roller length which touches the end of the wet area, while with  $Fr_{s_{\text{cr}}}=0.75$ a reduced roller is obtained.
    For both values of $Fr_{s_{\text{cr}}}$  the wave front is sharper than the one obtained with the other two criteria.
     As before, we compare the time evolution of the crest's free surface velocity with the wave's celerity. The results are summarized
in Figure \ref{fig_UC15Ph}. The breaking point of the local and hybrid criterion occur very early. The dissipation of the wave is such that
the curves never cross. This may be due to the effects of activating the limiter in the shock   smoothing of the waves peak.
Conversely,   the physical criterion provides a similar behaviour as the one observed in the previous case, with a blow-up of the particle velocity at the crest
after breaking onset. This is consistent with the sharper  wave fronts seen in Figure \ref{fig_Wei15break}.
  \vspace{-0.3cm}
\subsection{Breaking and Run-up of a Solitary Wave}
\label{Synolakis}
\begin{figure}[H]
\centering 
\includegraphics[width=0.48\textwidth]{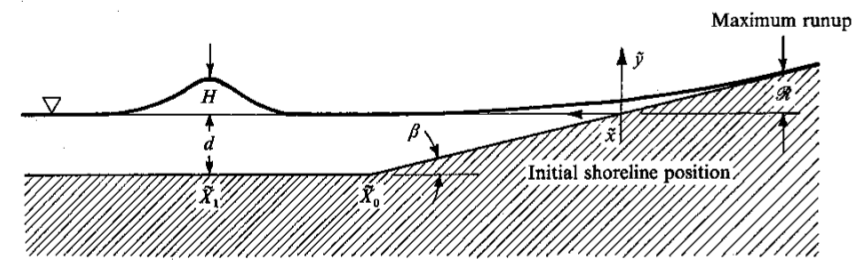}
\caption{Sketch of the test for the run-up of a solitary wave taken from \cite{Synolakis1987}}\label{sketch_syno}
\end{figure}

This is a classical benchmark due to by \cite{Synolakis1987}, which adds further complexity to the previous one, and considers the whole transformation of  a solitary  wave,  including the processes of breaking, run-up, and backwash.
The set-up consists in a solitary  wave  of  amplitude $A=0.28$m propagating on an  initial depth of $h_0=1$m and 
shoaling  on a 1:19.85  slope.  A sketch of the test is shown in Figure \ref{sketch_syno}.
\begin{figure}[H]
\centering
\subfigure[$t'=15$]{\includegraphics[width=0.49\textwidth]{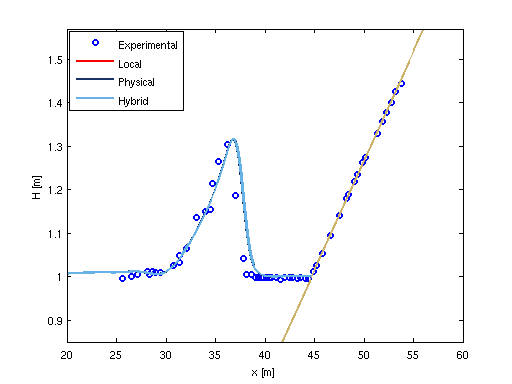}}
\subfigure[$t'=25$]{\includegraphics[width=0.49\textwidth]{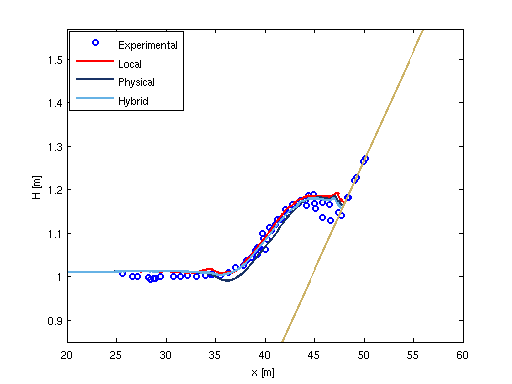}}\\
\subfigure[$t'=30$]{\includegraphics[width=0.49\textwidth]{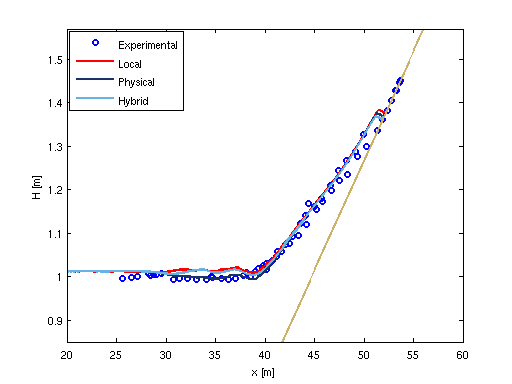}}
\subfigure[$t'=45$]{\includegraphics[width=0.49\textwidth]{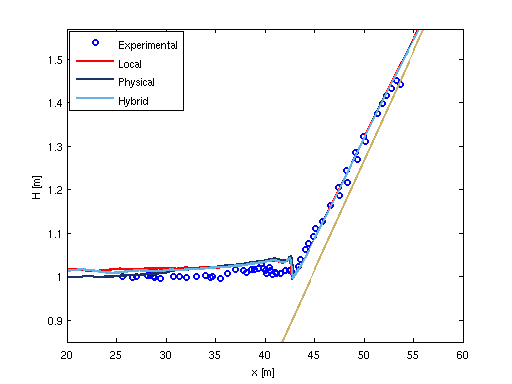}}\\
\subfigure[$t'=55$]{\includegraphics[width=0.49\textwidth]{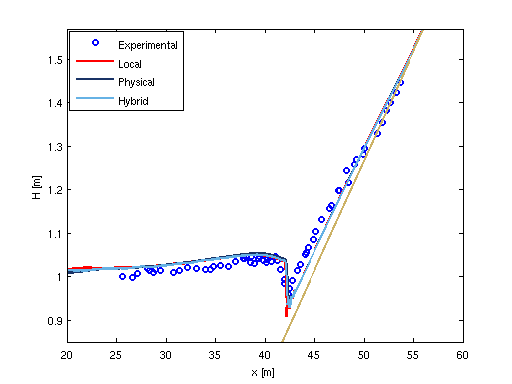}}
\subfigure[$t'=80$]{\includegraphics[width=0.49\textwidth]{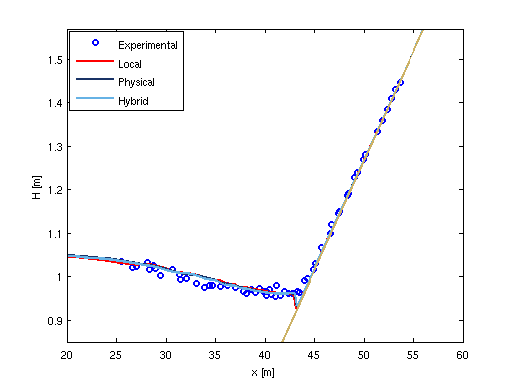}}
\caption{Results at different instants for the run-up of a solitary wave with $A=0.28$m for the hybrid model. Comparison between the local, hybrid, and physical breaking criteria. }\label{fig_Syno28_Bsqtot}
\end{figure}
Computations  have been run with a Manning coefficient of $0.01$, a CFL$=0.3$ and $\Delta x=0.1$m, i.e. the same parameters suggested in  \cite{Kazolea14}.
 The solution is visualized at several discrete instants, corresponding to   the beginning of the shoaling, the shoaling, the run-up and the back-wash. 
Results are displayed  w.r.t. the non-dimensional time  $t'=t\sqrt{g/h_0}$, and compared to experimental data.
The same problem has been used in section \S3 to evaluate the smoothness detection strategy proposed in the paper (see Figure \ref{fig_Syno_SW15}). 

The results are summarized in Figure  \ref{fig_Syno28_Bsqtot}. Despite the differences observed in the previous case, in this test  we can hardly see any distinction between the solutions 
obtained with the considered detection methods.  Overall, each approach shows an excellent agreement with the experimental data. 
We can see a slight phase advance at $t'=15$, which, however, is independent of the detection criterion, as none of the computations has started breaking at this point.
 This phase difference is assumed to be linked to the celerity of the exact solitary wave for \eqref{eq:Bsqsys_eq} given by (see  \cite{rf14} for details)
$$
c=\sqrt{gh_0}\sqrt{\Big(\frac{A}{h_0}\Big)^2\frac{1+\frac{A}{3h_0}}{\frac{A}{h_0}-\ln(1+\frac{A}{h_0})}} > \sqrt{gh_0}.
$$ 
The effect of this phase advance is not visible in the solutions at later times.  We underline the nice capturing of the hydraulic jump at $t'=45$, which validates 
the numerical limiting procedure proposed (cf. section \S3). All of the breaking criteria manage to capture this feature. The exception to this is the solution
obtained with the local criterion, that outlines some weak oscillations, revealing a weak flagging of this region.
Last, we would like to remark that for this type of problems a singularity is  encountered for the local criterion as   $h$ goes to zero when the slope crosses
the line corresponding to the initial/still water depth. This singularity is handled by the addition of a small quantity at the denominator, as discussed in sections \S4.1 and \S4.2.
More importantly, the value of $h$ in vicinity and after  this point is not relevant any longer for the definition of the physical quantities used in the breaking citeria.
In particular, we can note that the local criterion always breaks after this point, as  $|\eta|= |H-h| >|h|$ when $h<0$.
Similarly, the  first condition in the pre-flagging \eqref{secondcritI} loses its physical meaning, as $\sqrt{gh}$ no longer represents a relevant celerity.
For cases as the one in this section, this does not seem to affect the results. This is due to the fact that the three criteria end up detecting breaking in the run-up phase, and close to the 
backwash region where the depth is very small. These features are very well reproduced   by the Shallow Water equations.
However, as we will see later, for periodic waves, and especially in presence of strong set-up, this may be a source of inaccuracies.


\subsection{Periodic Wave over a Submerged Bar}
\label{Beji}  
In this benchmark  we  consider monochromatic  waves propagating over a submerged bar.  
According to \cite{Mingham2009}, this test originates from Dingemans \cite{dingemans1997} who verified Delft Hydraulics numerical model HISWA, and was then repeated by Beji and Battjes later \cite{Beji}.
\begin{figure}[H]
\centering
\includegraphics[scale=0.5]{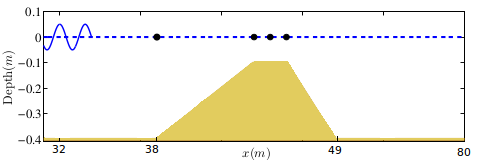}
\caption{Sketch of the submerged bar test case set-up.}
\label{fig_submergedBar}
\end{figure}
The test consists in a periodic wave with amplitude $A=0.027$m and period $T=2.525$s. The wave is  generated at $X_0=22$m over an initial water depth of $h_0=0.4$m.
The  bathymetry has  a  fore slope of 1:20, and a rear one of 1:10. The length of  the plateau is of $2$m, and the corresponding depth of $0.1$m.
Experimental data in four gauges (cf.  Figure \ref{fig_submergedBar}) are available for comparison. 
The problem is solved on the    domain   $[0, 60]$m, with  $\Delta x= 0.04$m, and absorbing conditions   on both ends. Periodic waves are generated using the source method of  \cite{rf14}.
This is   a difficult case, and  following \cite{Kazolea14}, we modify the values of  section \S4 for the detection criteria. We
 set $\mathcal{E}_{\text{critical}}=0.6$ for the local criterion, and  $\gamma=0.3$, $\phi=30\degree$  for the hybrid and physical ones.
For the latter, we will  use only  $Fr_{s_{\text{cr}}}=0.75$.

\begin{figure}[H]
\centering
\subfigure[ Gauge  $\text{G}_1$ - beginning of the upward slope]{\includegraphics[scale=.55]{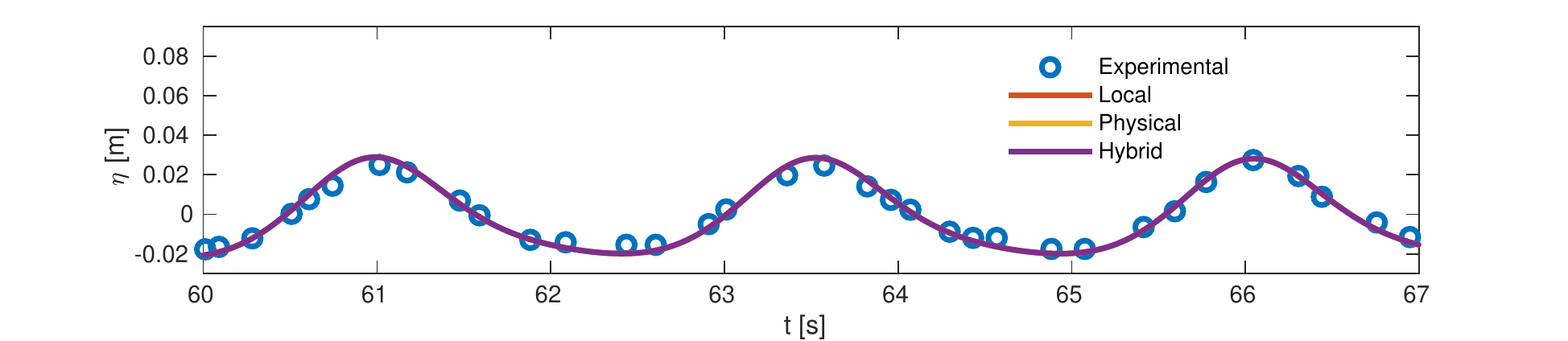}\label{fig_SubG1}}
\subfigure[ Gauge  $\text{G}_2$ - beginning of the plateau]{\includegraphics[scale=.55]{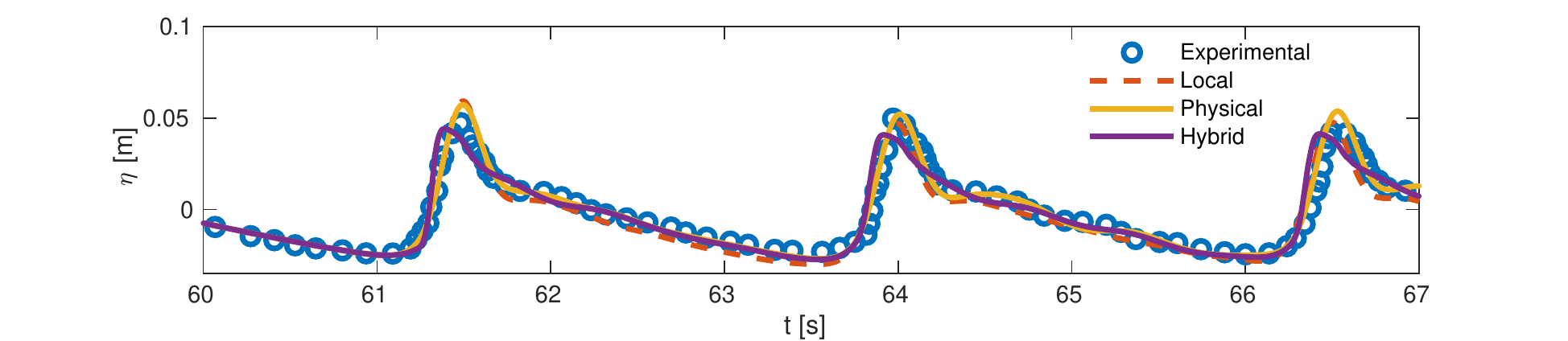}\label{fig_SubG2}}
\subfigure[Gauge 
$\text{G}_3$ - midpoint of the plateau]{\includegraphics[scale=.55]{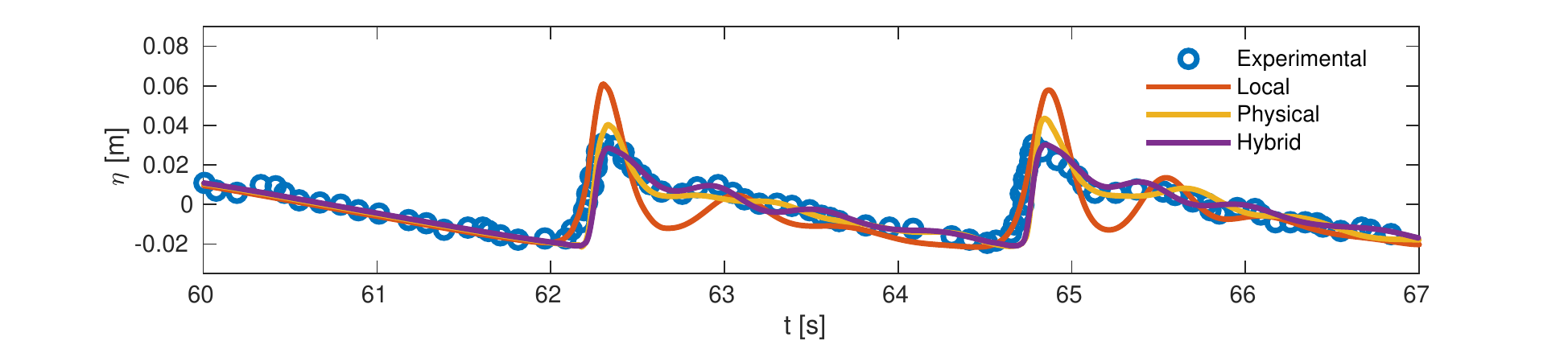}\label{fig_SubG3}}
\subfigure[ Gauge  $\text{G}_4$ - beginning of the downward slope]{\includegraphics[scale=.55]{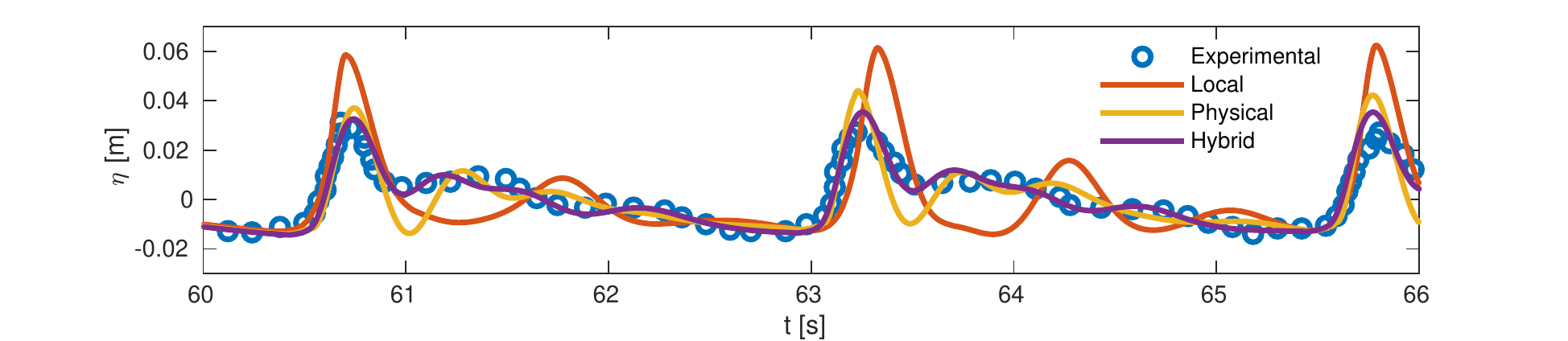}\label{fig_SubG4}}
\caption{Periodic waves over a  bar. Time series of the free surface elevation at gauges.}\label{fig_SubBARgauge}
\end{figure}

 For  $Fr_{s_{\text{cr}}}=1$   no breaking is observed at all.
Moreover,    reflections due to the boundary conditions may induce a  phase shift in the results, affecting the comparison with the experimental data \cite{KazoleaPM}.
For this reason, as in \cite{Kazolea14} we have focused our attention on the   first  generated waves,   before any interaction with the boundary conditions on the right end of the domain 
has taken place. In practice,  the first gauge is used to re-phase the numerical signals, and the resulting phase calibration is applied to all gauges.
The resulting water elevation time series in the three remaining gauges are reported in Figure  \ref{fig_SubBARgauge}. 

The  comparison shows that the hybrid criterion is the only one capturing the correct magnitude and phase of the waves.
The local criterion essentially fails to capture the onset of breaking, even with the reduced value of $\mathcal{E}_{\text{critical}}=0.6$. 
This leads to wrong amplitudes as well as to a shift in the phase of  parts of the signal. The physical criterion gives acceptable results, managing  to capture the phase of the experimental
data. The amplitudes remain slightly higher than those of the data and obtained with the hybrid detection. This may be the consequence of a later onset
of breaking,  as well as of some intermittency  in the breaking process, reducing the overall dissipation. 

\begin{figure}[H]
\centering
\subfigure[Local ]{\hspace{-0.9cm}\includegraphics[width=0.4\textwidth]{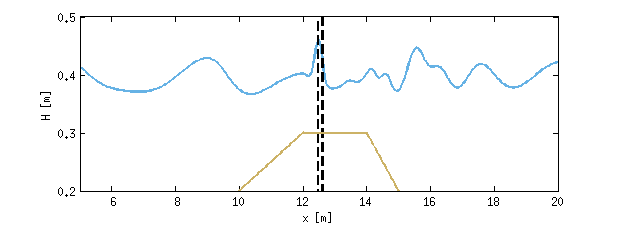}
\hspace{-0.5cm}\includegraphics[width=0.4\textwidth]{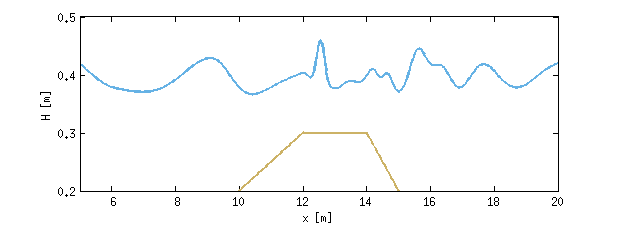}
\hspace{-0.5cm}\includegraphics[width=0.4\textwidth]{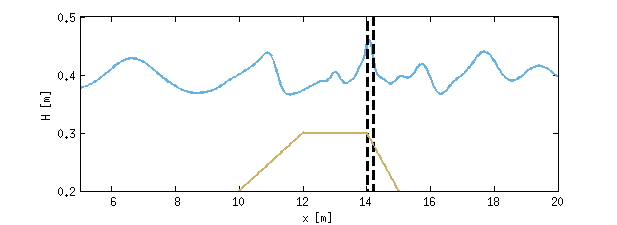}\label{fig_SubTP}}
\subfigure[Hybrid ]{\hspace{-0.9cm}\includegraphics[width=0.4\textwidth]{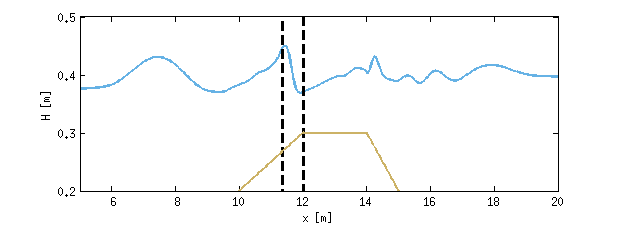}
\hspace{-0.5cm}\includegraphics[width=0.4\textwidth]{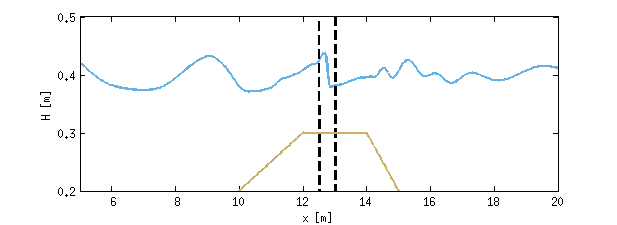}
\hspace{-0.5cm}\includegraphics[width=0.4\textwidth]{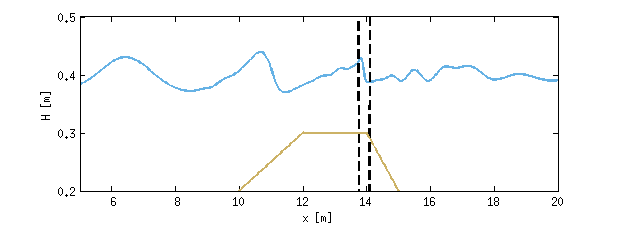}\label{fig_SubDelis}}
\subfigure[Physical  ]{\hspace{-0.9cm}\includegraphics[width=0.4\textwidth]{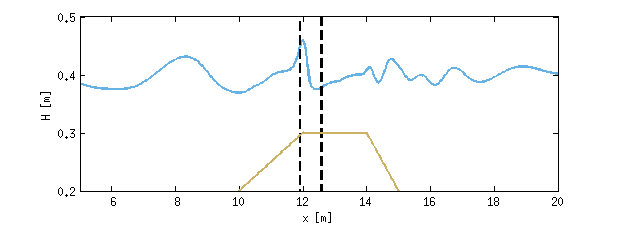}
\hspace{-0.5cm}\includegraphics[width=0.4\textwidth]{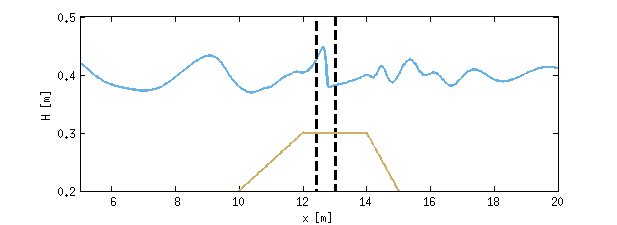}
\hspace{-0.5cm}\includegraphics[width=0.4\textwidth]{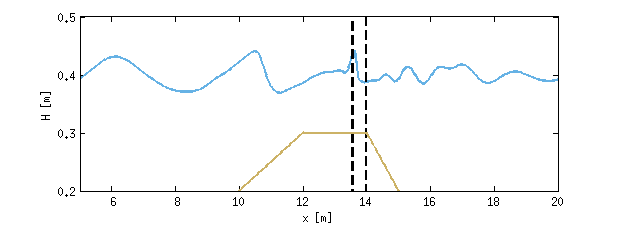}\label{fig_SubPh}}
\caption{Periodic wave over a submerged bar. Comparison of the instances when the breaking first occurs on a wave (\text{left}) at a mid-time (\text{center}) and the last breaking point (\text{right}).}\label{fig_SubBARcomp}
\end{figure}
To better understand the differences in the time series, we have plotted the wave profiles  at  time instants corresponding to the first and last breaking instants on the plateau.
Note that these differ for the three criteria. We also have included the wave profile for all the solutions at  the same intermediate time.
The resulting free surface distributions are displayed in Figure  \ref{fig_SubBARcomp}. 

If we set $t=0$  in correspondence of  the first breaking  of the hybrid criterion,
the physical criterion breaks roughly at  $0.18\,T$, while the local detection gives an onset at 
$t=0.345\,T$, with $T$ the incoming wave period. The delay is clearly visible in the left pictures in Figure \ref{fig_SubBARcomp}. 
The plots also show that the local criterion still gives a very small roller, while the hybrid and physical criterion gives a breaking region of comparable size.
The central pictures, referring to $t= 0.365\,T$, show that the local criterion has already deactivated the flagging, 
while  both the hybrid and physical criteria are still breaking. Before  the wave leaves the plateau, 
the local criterion activates again at  $t\approx 0.82 \,T$, probably due to the decrease in bathymetry. This sudden flash in breaking has  no effect on the wave shape. 
The hybrid and the physical criteria keep the flagging activated until $t\approx 0.77 \,T$s and $t\approx 0.70 \,T$, respectively.   
One must however also consider that in general, at least in our implementation, the  lack of memory in the flagging introduces a degree of intermittency.
This is observed  more  for the physical criterion than for the hybrid one. The cause may be the noise in the evaluation of the derivatives defining $u_s$ (cf, section \S 4.3),
and the result is a weaker dissipation leading to the slightly higher amplitudes of Figure   \ref{fig_SubBARgauge}.

\subsubsection*{\revP{Mesh convergence study}}
\revP{We use this test to evaluate the mesh dependence of the closure and detection criteria. Note that a very thorough investigation of this issue has already been presented in \cite{kazolea2018}.
The reference, as well as our study shows that the type of breaking closure based on the coupling of Boussinesq-type equations with the Non-linear Shallow Water model and the approximation of the rollers by bores,  introduces instabilities, that lead to a blow up, as the mesh gets refined. This behavior corresponds most likely to a mismatch between the energy transported in the dispersive field which is not entirely transmitted at the interface. Smoothing the coupling numerically only delays this blow up.
Note that similar results, albeit with different schemes, are reported in \cite{shi2012,filippini2016,kazolea2018}. }

\revP{ The grid convergence study reported in this section allows to include in the grid convergence analysis the effect of using different detection criteria.
For the scheme and model used here, grid convergence for smooth, non-breaking, cases can be already found in \cite{rf14,filippini2016}.
We thus repeat the computations of the submerged bar by halving the mesh size from $\Delta x=0.04$m, used in the previous section, to $\Delta x=0.02$m and $\Delta x=0.01$m.
Only, the local criterion has allowed to compute the whole solution with an additional refinement of $\Delta x=0.005$m. }

\begin{figure}[h!]
\centering
\subfigure[Local Criterion]{\includegraphics[scale=.55]{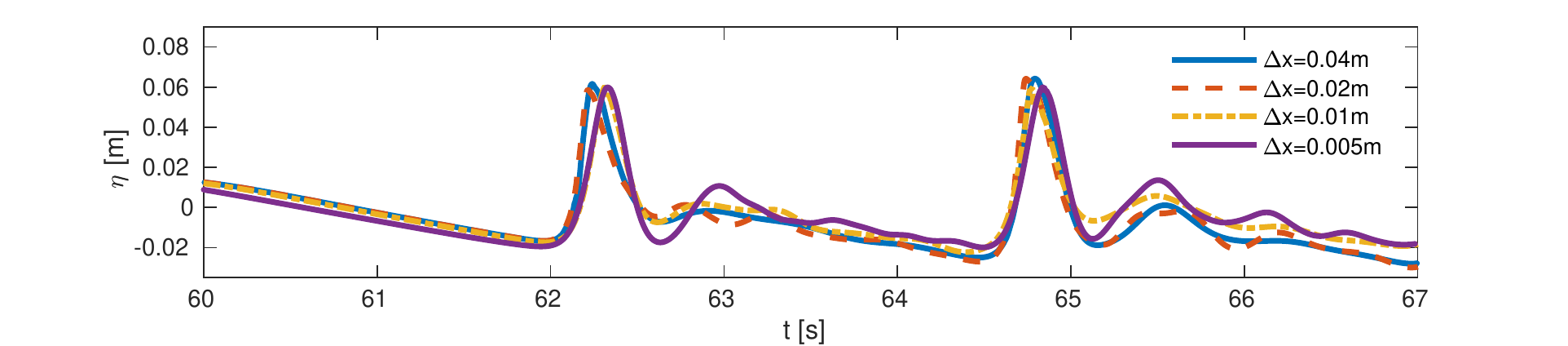}\label{fig_SubG3_TP}}
\subfigure[Physical Criterion]{\includegraphics[scale=.55]{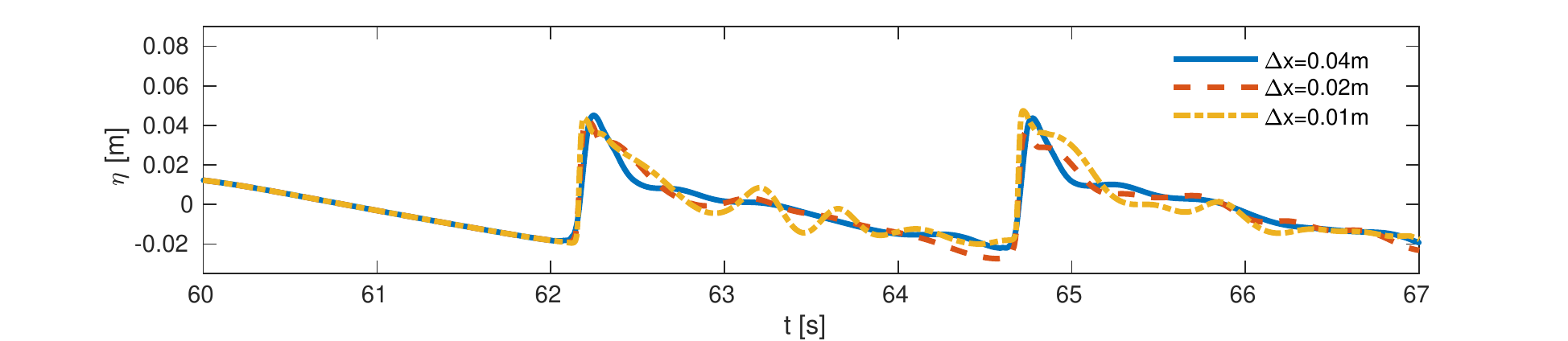}\label{fig_SubG3_FR}}
\subfigure[Hybrid Criterion]{\includegraphics[scale=.55]{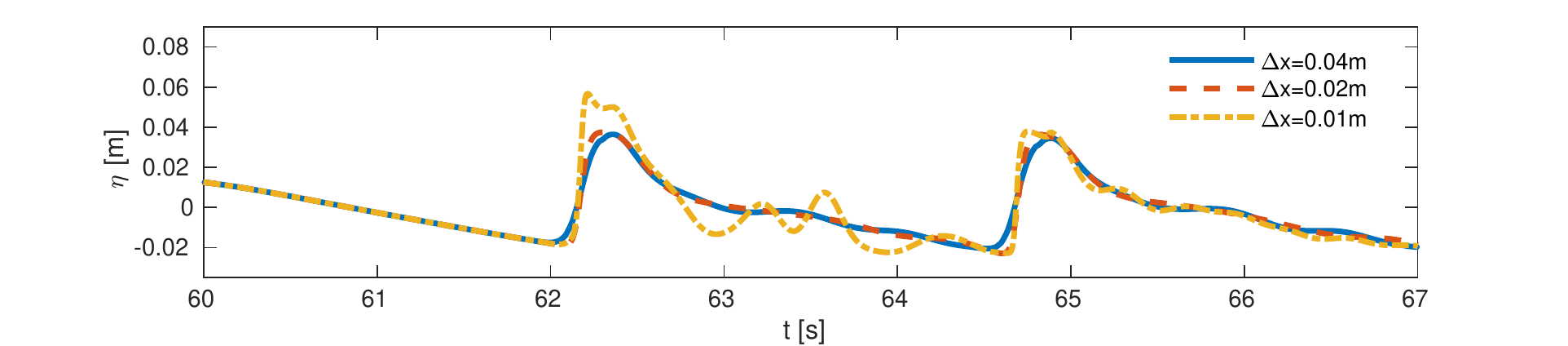}\label{fig_SubG3_PH}}
\caption{\revP{Periodic waves over a bar. Time series of the free surface elevation at gauge $G_3$ for each criteria with different mesh sizes.}}\label{fig_SubBARgauge_meshconv}
\end{figure} 

\begin{figure}[h!]
\centering
\includegraphics[scale=.55]{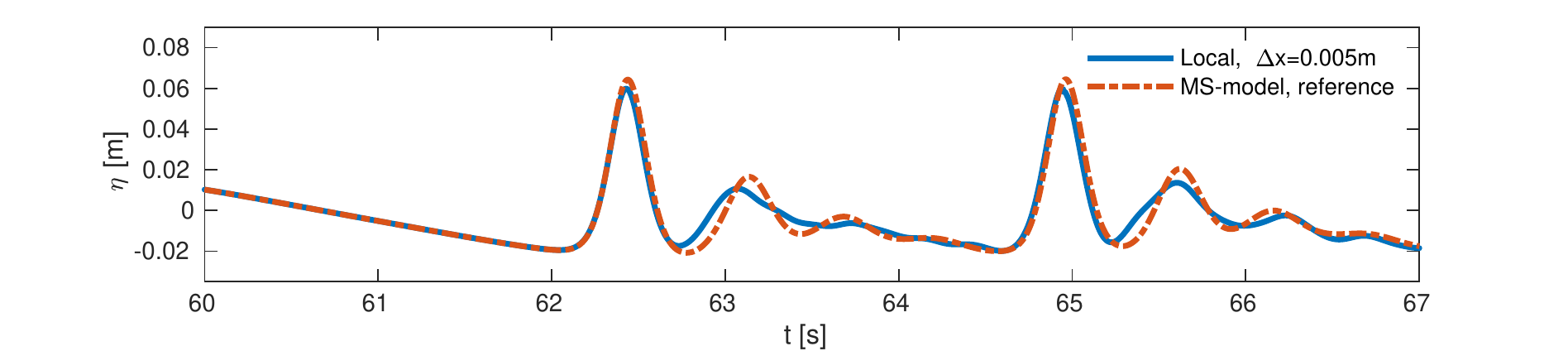}
\caption{\revP{Periodic waves over a bar. Comparison on gauge $G_3$. Local criterion on $\Delta x=0.005$m vs.  pure Madsen-S{\o}rensen (MS) model on $\Delta x=0.01$m.}}\label{fig_SubBARgauge_tp_bsq}
\end{figure}

\revP{On Figure \ref{fig_SubBARgauge_meshconv} we report the results for the water levels in gauge $G_3$, but the trends observed in the other gauges are similar.
In the figure, the first row reports the results obtained with the local criterion, the second those of the physical detection criterion, and the last the results obtained with the 
hybrid criterion. For the last two criteria our results are very similar to those of \cite{kazolea2018}, and show the appearance on the finer meshes of unphysical high frequency
oscillations superimposed to the trend of the coarse mesh solution. This phenomenon is a not as visible for the results obtained with the local criterion, which only flashes at the end of the upward slope. }
\revP{Compared to the weakly nonlinear model analyzed in \cite{kazolea2018} the oscillations obtained here seem slightly weaker in both in frequency and amplitude.
We believe this to be a consequence of the shoaling properties of the Madsen-S{\o}rensen equations used here which have a tendency to under-shoal, compared to the Nwogu-model
of the reference, and to give lower and longer waves, especially in the nonlinear range \cite{fbcr15}. Concerning the detection criteria, the main conclusions we can draw are that
the onset and growth of the instabilities seem to be more visible with the hybrid detection which for this case is the one breaking more strongly. The phenomenon is still very 
intense with the physical criterion, which still detects wave breaking on most of the plateau, and much less present with the local detection criterion which breaks very weakly.
Indeed, the local criterion, for the mesh $\Delta x=0.005$m appears to have the same shape as a Boussinesq-type solution, i.e. without any breaking applied at all, as can be seen from Figure \ref{fig_SubBARgauge_tp_bsq}.}

\begin{figure}[h!]
\centering
\subfigure[Physical Criterion]{\includegraphics[scale=.55]{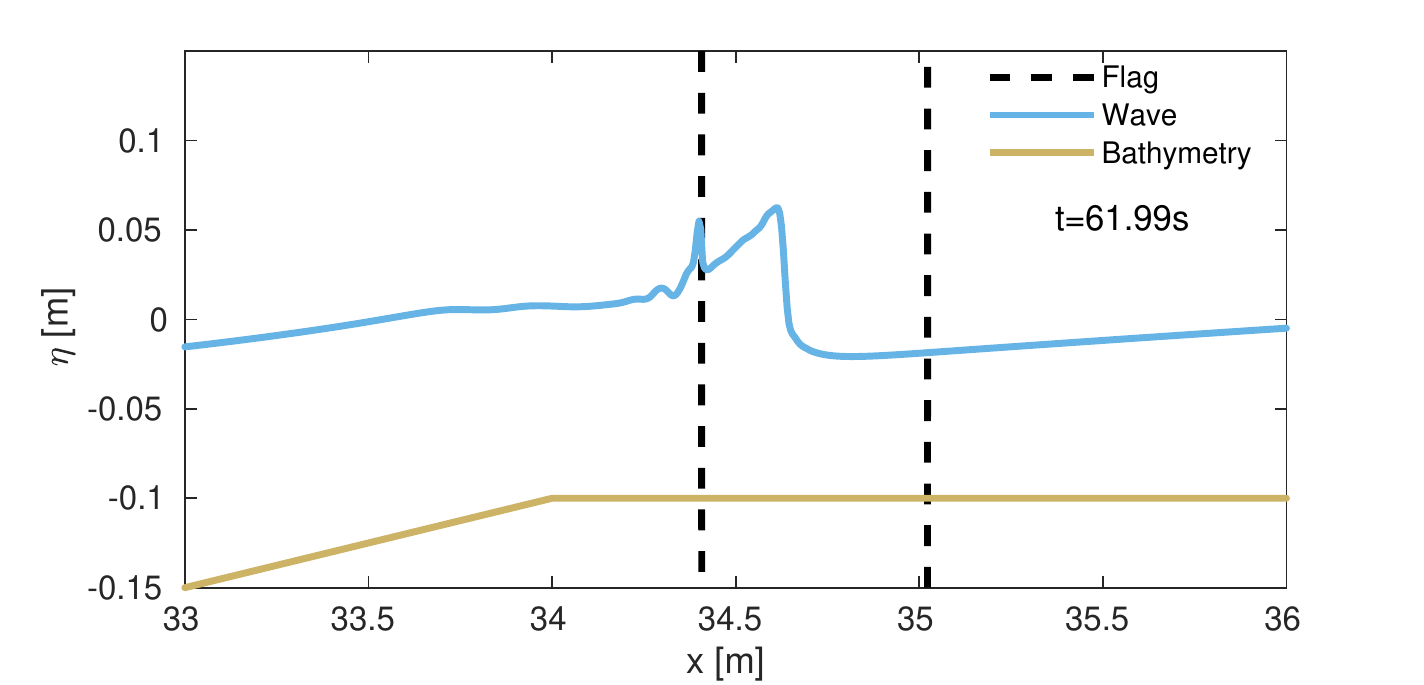}\label{fig_SubG3_FR1}}
\subfigure[Hybrid Criterion]{\includegraphics[scale=.55]{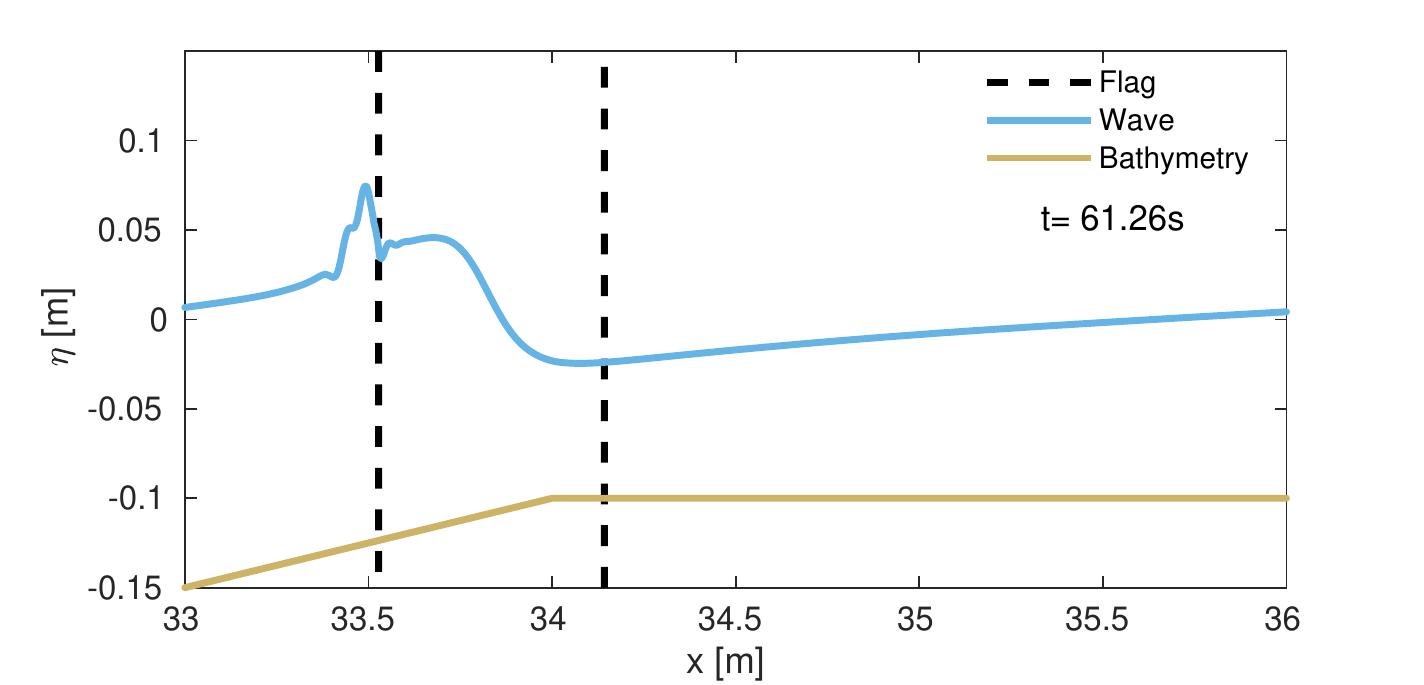}\label{fig_SubG3_PH1}}
\caption{\revP{Periodic waves over a bar. Snapshots with first high-frequency oscillations at $\Delta x=0.005$m.}}\label{fig_SubBARgauge_fail}
\end{figure}

\revP{Finally, as a last comparison, in Figure \ref{fig_SubBARgauge_fail} we show how the physical and hybrid criterion start failing with the very first breaking instances, when the mesh gets too fine, as with here $\Delta x=0.005$m. }

\subsection{Breaking and set-up of Monochromatic Waves on a Shore}
\label{Hansen}
As a last benchmark, we consider the shoaling and breaking  of monochromatic waves over a shore, based on the experiments of 
 Hansen and Svendsen \cite{hansen79}.  This is an interesting case, as the set-up is induced by the breaking near the shore.
The waves considered are generated on a still water depth  $h_0=0.36$m and then propagate on a slope of 1:32.26. 
As in \cite{Kazolea14}  the   grid size  used is $\Delta x=0.025$m.   Wave generation, as well as the  sponge layers  on the left  end of the domain, are set up following \cite{rf14}.

\begin{figure}[H]
\centering
\subfigure[Local  ]{\hspace{-0.8cm}\includegraphics[width=0.5\textwidth]{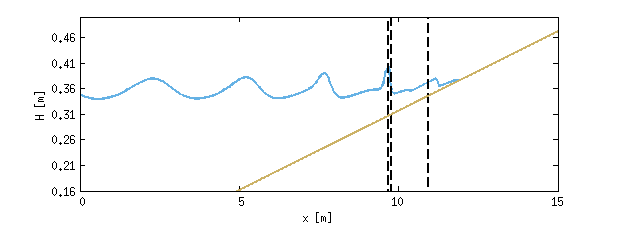}
\hspace{-0.725cm}\includegraphics[width=0.5\textwidth]{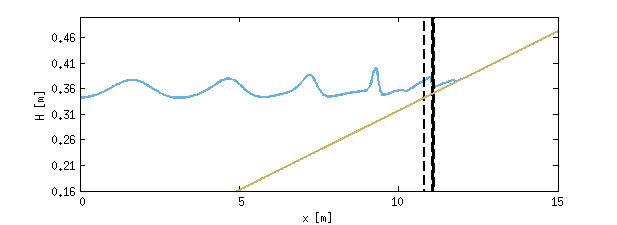}\label{HSlocal}}
\subfigure[Hybrid  ]{\hspace{-0.8cm}\includegraphics[width=0.5\textwidth]{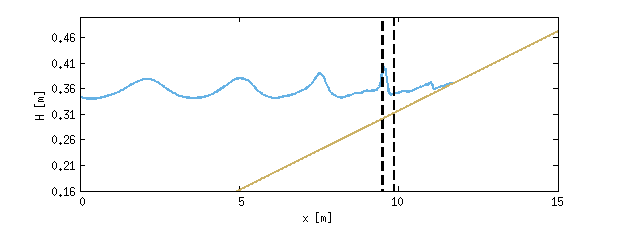}
\hspace{-0.725cm}\includegraphics[width=0.5\textwidth]{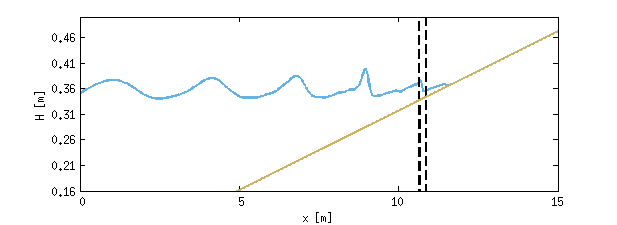}\label{HShybrid}}
\subfigure[Physical   ($Fr_{s_{\text{cr}}}=1$) ]{\hspace{-0.8cm}\includegraphics[width=0.5\textwidth]{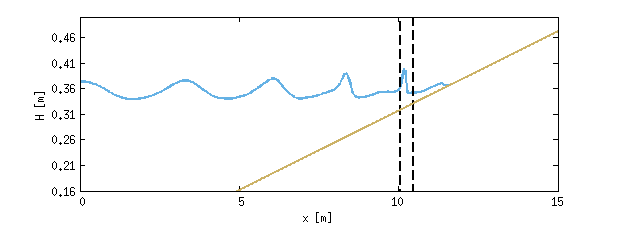}
\hspace{-0.725cm}\includegraphics[width=0.52\textwidth]{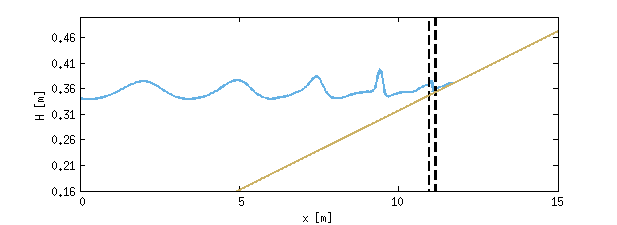}\label{HSphysical}}
\subfigure[Physical   ($Fr_{s_{\text{cr}}}=0.75$) ]{\hspace{-0.8cm}\includegraphics[width=0.5\textwidth]{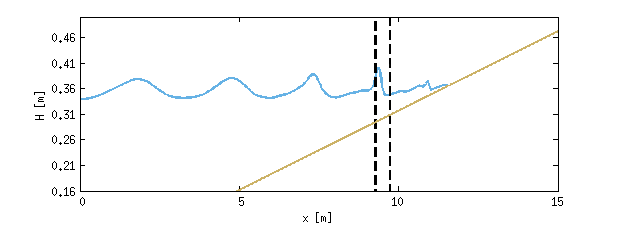}
\hspace{-0.725cm}\includegraphics[width=0.5\textwidth]{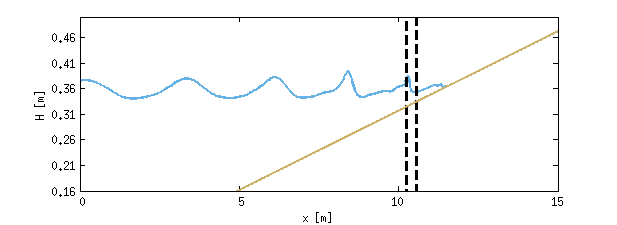}\label{HSphysical75}}
\caption{Hansen and Svendsen test 051041. Snapshots of the first and last breaking times.}\label{HansenSorensen}
\end{figure}

{\bf Case 051041.} For this case waves have  a period $T=2$s and amplitude  $A=0.018$m, corresponding to an Ursell number  $U=4.8077$.
We visualize in Figure \ref{HansenSorensen} the wave profiles at the   first and last breaking instants (for a fixed wave) detected by each criteria. 

\begin{figure}[H]
\centering
\subfigure[Local ]{\includegraphics[width=0.5\textwidth]{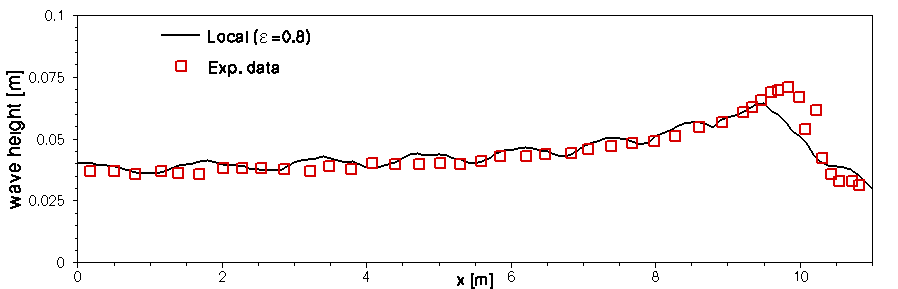}
\includegraphics[width=0.5\textwidth]{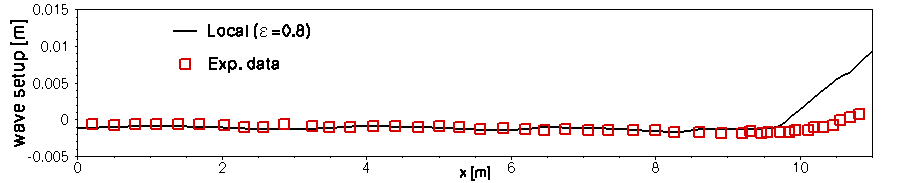}\label{HS1_tp}}
\subfigure[Hybrid ]{\includegraphics[width=0.5\textwidth]{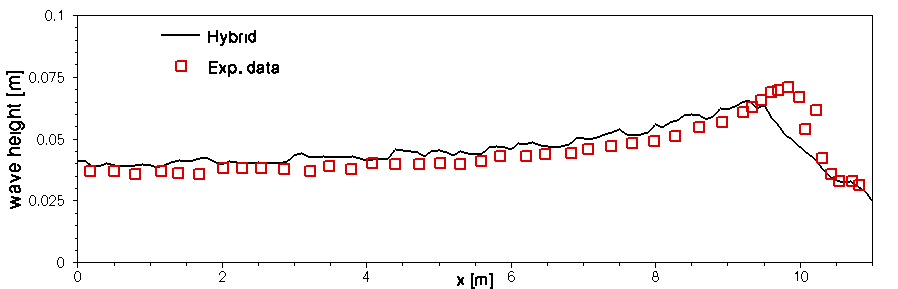}
\includegraphics[width=0.5\textwidth]{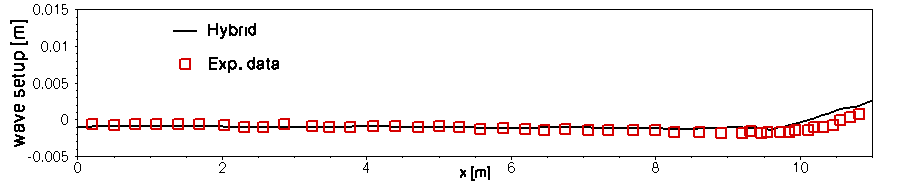}\label{HS1_D}}
\subfigure[Physical  ($Fr_{s_{\text{cr}}}=1$)]{\includegraphics[width=0.5\textwidth]{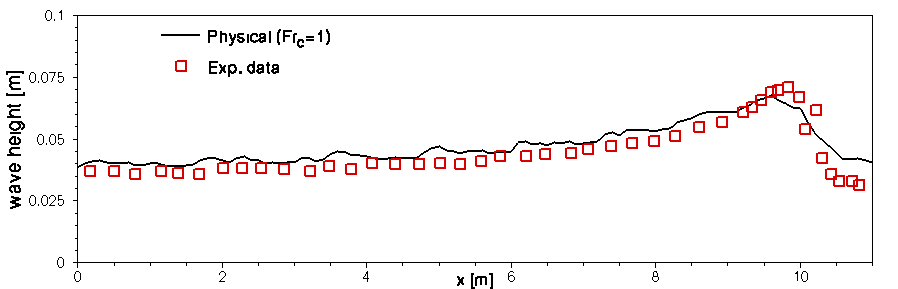}
\includegraphics[width=0.5\textwidth]{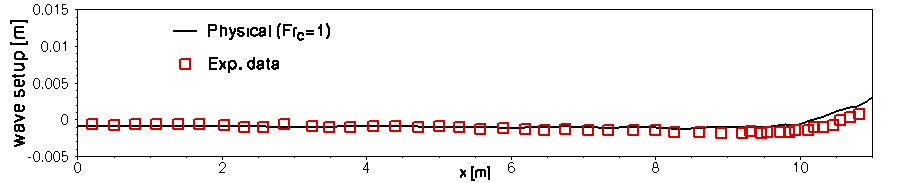}\label{HS1_Ph1}}
\subfigure[Physical  ($Fr_{s_{\text{cr}}}=0.75$)]{\includegraphics[width=0.5\textwidth]{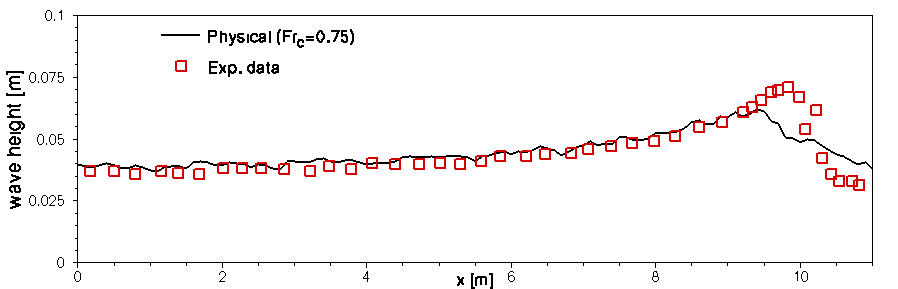}
\includegraphics[width=0.5\textwidth]{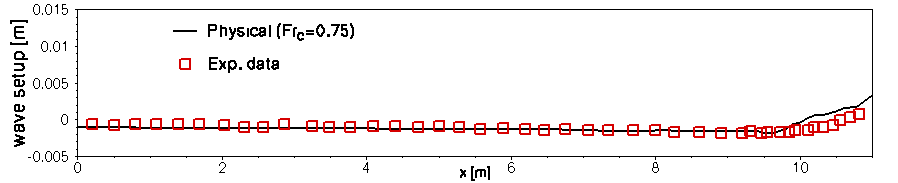}\label{HS1_Ph75}}
\caption{Hansen and Svendsen - 051041. Wave height \textit{(left)} and mean water level \textit{(right)}.}\label{HansenSorensensetup}
\end{figure}
As before, the breaking region is highlighted by vertical lines.  As in the other cases, the onset of breaking is predicted slightly  later by the physical criterion 
when $Fr_{s_{\text{cr}}}=1$. 
Differently than in the other benchmarks, here for $Fr_{s_{\text{cr}}}=0.75$ the onset is predicted earlier than with all other methods.
The termination is also shifted accordingly for the two cases.
We can also see again that the roller size computed by the local criterion is much shorter. For the latter, we can see a single vertical line at $x\approx 11.5$, corresponding to the zero of $h$,
where breaking is turned on until the end of the domain. As discussed at the end of section \S5.3  this is a consequence of the use of $h$ in the  definition of a reference depth in the breaking detection.
In Figure  \ref{HansenSorensensetup}, we compare the spatial distribution of the wave height and the mean water level (set-up)  with the experimental data of \cite{hansen79}.
The results are qualitatively close to those reported in \cite{Kazolea14}.  The prediction of the wave heights is clearly  affected from the prediction of the breaking onset.
Even though under-predicting  the breaking wave height, the physical criterion with $Fr_{s_{\text{cr}}}=1$ gives 
the best result for this case  in terms of position of the breaking onset and slope of the height after the  breaking point. 
All the other criteria, including the physical one for $Fr_{s_{\text{cr}}}=0.75$, provide a very early onset of breaking, with an under-prediction of the breaking wave height, 
and a milder slope after.  Some of the height underprediction may be also related to the under-shoaling characteristics of the Madsen and S{\o}rensen model \eqref{eq:Bsqsys_eq}.
Looking at the wave set-up,   reported on the right in the same figure, we can clearly see the effects of the use of $h$ in the local criterion. An unrealistically high set-up is observed. As shown in \cite{Bonneton2007}  this can be   related to 
the amount of dissipation injected  in the region where $h<0$ which is entirely flagged as breaking. 
This has no effects on the results obtained with the  other criteria, providing a correct set-up.

\begin{figure}[H]
\centering
\subfigure[Local ]{\includegraphics[width=0.5\textwidth]{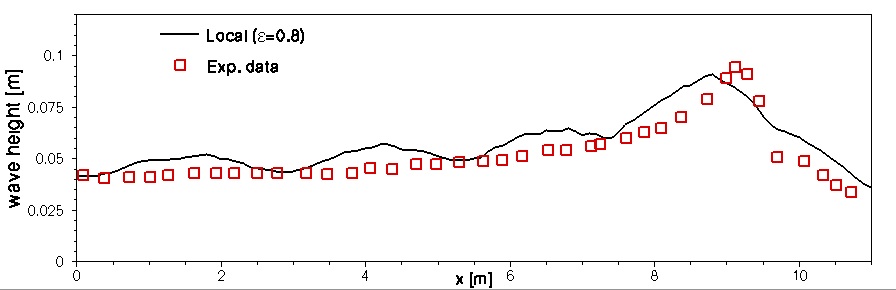}
\includegraphics[width=0.5\textwidth]{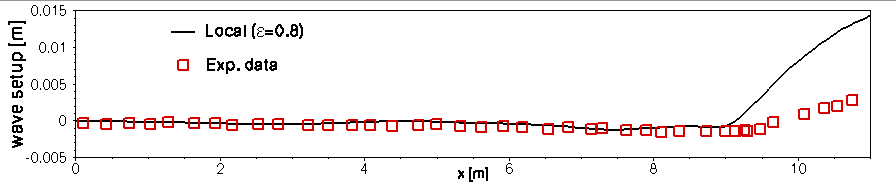}\label{HS2_tp}}
\subfigure[Hybrid ]{\includegraphics[width=0.5\textwidth]{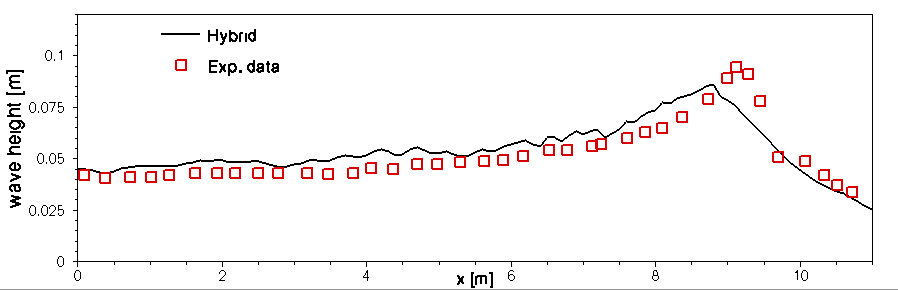}
\includegraphics[width=0.5\textwidth]{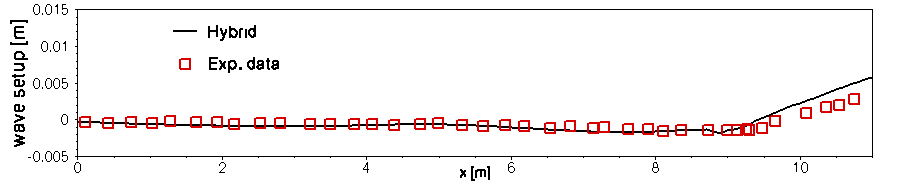}\label{HS2_D}}
\subfigure[Physical  ($Fr_{s_{\text{cr}}}=1$)]{\includegraphics[width=0.5\textwidth]{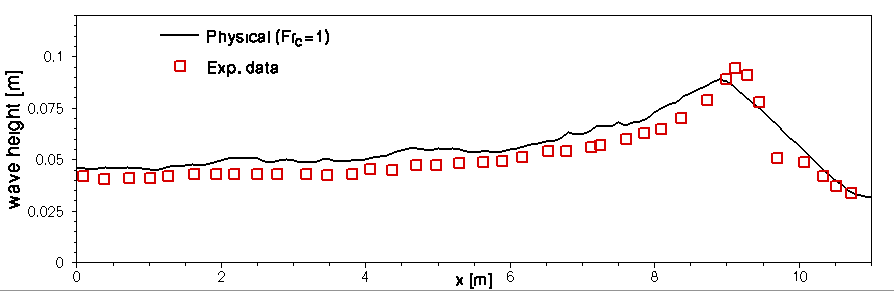}
\includegraphics[width=0.5\textwidth]{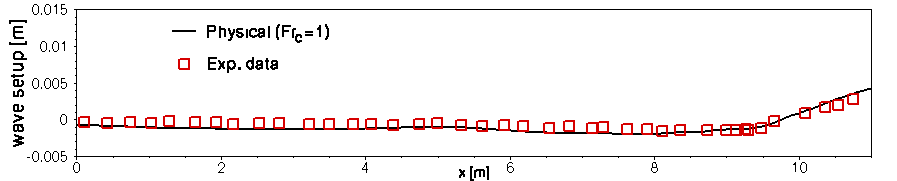}\label{HS2_Ph1}}
\subfigure[Physical  ($Fr_{s_{\text{cr}}}=0.75$)]{\includegraphics[width=0.5\textwidth]{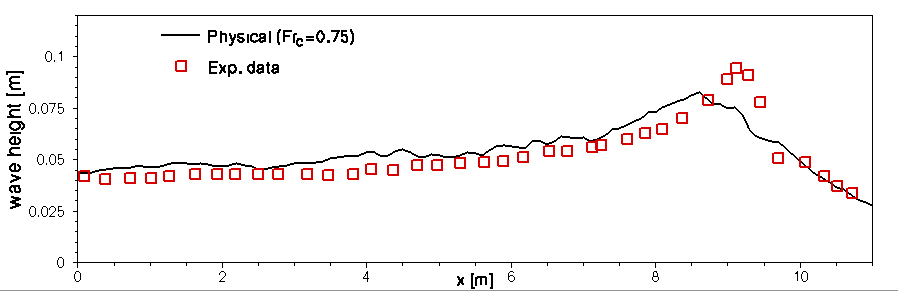}
\includegraphics[width=0.5\textwidth]{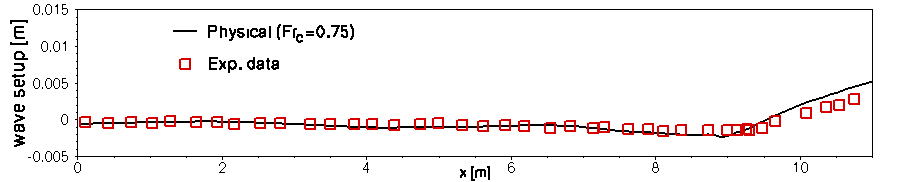}\label{HS2_Ph75}}
\caption{Hansen and Svendsen test 031041. Wave height \textit{(left)} and mean water level \textit{(right)}.}\label{HansenSorensensetup2}
\end{figure}

{\bf Case 031041.} In this case waves have a period $T=3.333333$s and amplitude  $A=0.0215$m, corresponding to an Ursell number    $U=17.5588$.
  \begin{figure}[H]
\centering
\subfigure[Local  ]{\hspace{-0.8cm}\includegraphics[width=0.5\textwidth]{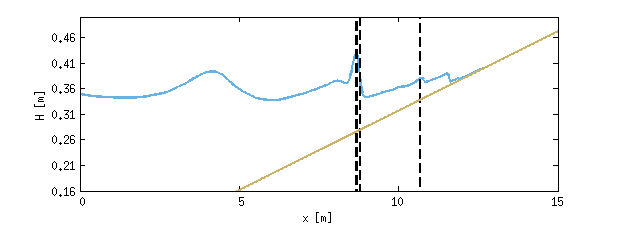}
\hspace{-0.725cm}\includegraphics[width=0.5\textwidth]{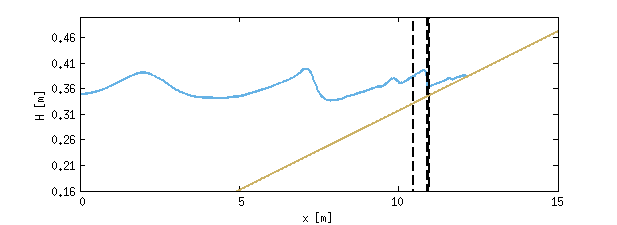}\label{HS2local}}
\subfigure[Hybrid  ]{\hspace{-0.8cm}\includegraphics[width=0.5\textwidth]{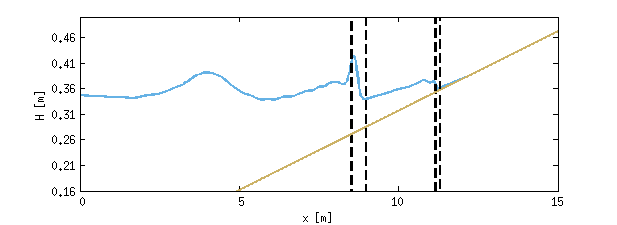}
\hspace{-0.725cm}\includegraphics[width=0.5\textwidth]{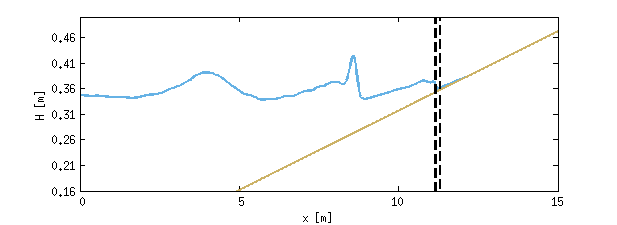}\label{HS2hybrid}}
\subfigure[Physical   ($Fr_{s_{\text{cr}}}=1$) ]{\hspace{-0.8cm}\includegraphics[width=0.5\textwidth]{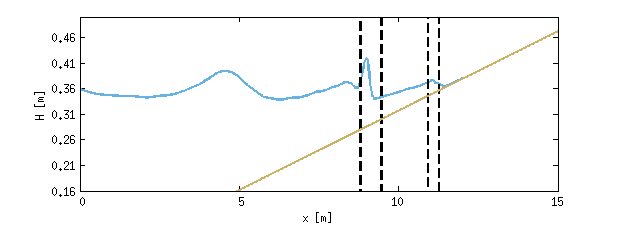}
\hspace{-0.725cm}\includegraphics[width=0.5\textwidth]{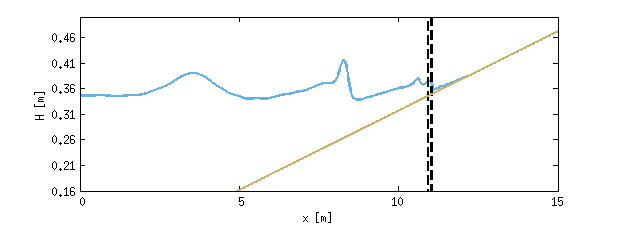}\label{HS2physical}}
\subfigure[Physical   ($Fr_{s_{\text{cr}}}=0.75$) ]{\hspace{-0.8cm}\includegraphics[width=0.5\textwidth]{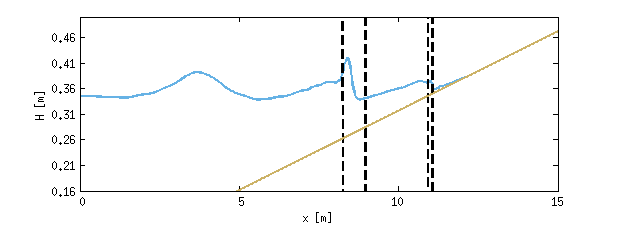}
\hspace{-0.725cm}\includegraphics[width=0.5\textwidth]{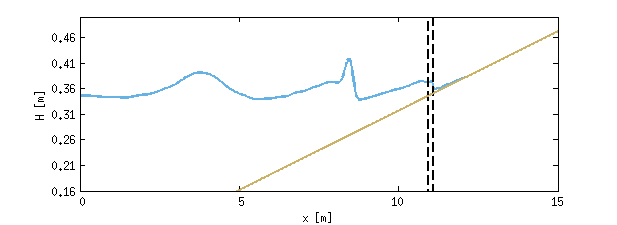}\label{HS2physical1}}
\caption{Hansen and Svendsen test 031041.  Snapshots of the first and last instant of breaking for different detection criteria.}\label{HansenSorensen2}
\end{figure}
The results  are summarized in Figures \ref{HansenSorensen2}, and   \ref{HansenSorensensetup2},
 reporting  respectively the  wave profiles at breaking onset/termination for a given wave, and the mean water levels and  set-up.
We can observe a behaviour similar to the one discussed in the previous tests. The local criterion gives very thin rollers, and flags as breaking occurs within the entire region where the bathymetry is above the zero of $h$. 
The onset is predicted slightly  later for the physical criterion if  $Fr_{s_{\text{cr}}}=1$, but the difference is less important in this case. Concerning the mean water level, the best is again the one
 obtained with   the physical criterion and  $Fr_{s_{\text{cr}}}=1$, the other methods, including the physical criterion for  a lower value of the critical Froude,  are  affected by the early prediction of breaking onset.

The same drawbacks of the local criterion can be observed in the wave set-up, which is too large for this detection method, at least with the implementation discussed in section \S4.
A good prediction of this quantity is instead provided by the hybrid and physical criteria.

\newpage
\section{Conclusions and Outlook} 
We have discussed the numerical implementation and the extensive benchmarking of three breaking detection criteria in a finite element Boussinesq model.
Breaking closure is handled by means of the  hybrid method, that reverts the enhanced Boussinesq equations to the hydrostatic Shallow Water system,  and models wave breaking fronts as  shocks.
The detection criteria tested involve an original implementation in a numerical scheme of a more physical detection based on the comparison of the free surface particle velocity at the peak of the wave, with the wave celerity. 
\revP{In the case of the convective criterion, the lowering of the $Fr_{s_{\text{cr}}}$ to $0.75$ has been set empirically, to account for the under-shoaling of the   Boussinesq model \cite{fbcr15}.}
The thorough benchmarking performed shows an advantage in using detection methods  based on more sound physical arguments. 
In its current implementation the physical criterion is comparable to the hybrid slope/vertical velocity method proposed in  \cite{Kazolea14},
which is already quite satisfactory in a wide range of cases. It also allows to obtain improved results  in some cases, as  the Hansen and Svendsen  experiments of section \S5.4. 
The results also show some limitations related to the adequacy of the detection criteria to  the propagation model.
More specifically, the classical asymptotic formula used for the vertical velocity profile may be less well adapted for the type of Boussinesq model used than for, e.g.,  such as the   standard Peregrine equations of the Beji and Nadaoka \cite{fbcr15}. 
This is clearly an aspect which could be refined and investigated in the future, also including 
fully non-linear models.
\revP{Further, the switch between the detection criteria undergoes some smoothing but this will not guarantee the stability of the coupling as $\Delta x \rightarrow 0$.
This has been already investigated in \cite{filippini2016,kazolea2018} for the Green-Naghdi and Nwogu Boussinesq equations. 
The grid convergence study performed here confirms this observation also for the Madsen and Sorensen Model. The study also allows
to  give some highlight on the impact of the nature of flagging on this aspect. In particular  more than the intermittency of the flagging, the results
suggest that  the stronger the duration of the  breaking process the earlier the breakdown.  When convergence is obtained, as in the case 
of the local criterion in our test, the converged solution is essentially the non-breaking Boussinesq one, which we have shown quantitatively. 
The breakdown of the solution as the mesh is refined is 
a flaw for this closure, which deserves  a specific investigation.}\\
 Another element which deserves further study is the definition of the front celerity, which becomes less trivial in multiple space dimensions.
This may require to combine the ideas explored here with techniques used in image processing, as well as, in wind waves analysis. The main challenge will be to 
retain the locality of the detection, thus avoiding non-local operations as, e.g., Fourier transforms. 
Other elements more specifically linked to the numerical method also deserve some further refinement.
It is the case for example for the approximation of the second derivatives necessary in the vertical velocity development.
These may require special attention in the multidimensional case, especially on unstructured grids. 
The use of methods based on high-order polynomial expansions (see \cite{ar17}) may help, but some form of filtering will certainly  still be necessary.


\bibliographystyle{spmpsci}      
\bibliography{biblio}

\begin{thebibliography}{10}
\providecommand{\url}[1]{{#1}}
\providecommand{\urlprefix}{URL }
\expandafter\ifx\csname urlstyle\endcsname\relax
  \providecommand{\doi}[1]{DOI~\discretionary{}{}{}#1}\else
  \providecommand{\doi}{DOI~\discretionary{}{}{}\begingroup
  \urlstyle{rm}\Url}\fi

\bibitem{ar17}
Abgrall, R., Ricchiuto, M.: {High-Order Methods for CFD}.
\newblock pp. 1--54. Wiley Online Library (2017)

\bibitem{bacigaluppi:hal-01087945}
Bacigaluppi, P., Ricchiuto, M., Bonneton, P.: {A 1D Stabilized Finite Element
  Model for Non-hydrostatic Wave Breaking and Run-up }.
\newblock In: J.~Fuhrmann, M.~Ohlberger, C.~Rohde (eds.) {Finite Volumes for
  Complex Applications VII-Elliptic, Parabolic and Hyperbolic Problems}, pp.
  779--790. Springer International Publishing, Cham (2014)

\bibitem{bacigaluppi:hal-00990002}
Bacigaluppi, P., Ricchiuto, M., Bonneton, P.: {Upwind Stabilized Finite Element
  Modelling of Non-hydrostatic Wave Breaking and Run-up}.
\newblock Research Report RR-8536, {INRIA} (2014)

\bibitem{Beji}
Beji, S., Battjes, J.: Experimental investigation of wave propagation over a
  bar.
\newblock Coastal Engineering \textbf{19}(1-2), 151--162 (1993)

\bibitem{Bjorkavag2011}
Bj{\o}rkav{\aa}g, M., Kalisch, H.: Wave breaking in {B}oussinesq models for
  undular bores.
\newblock Physics Letters A \textbf{375}(14), 1570--1578 (2011)

\bibitem{Bonneton2007}
Bonneton, P.: Modelling of periodic wave transformation in the inner surf zone.
\newblock Ocean Engineering \textbf{34}(10), 1459--1471 (2007)

\bibitem{bfwts06}
Borthwick, A., Ford, M., Weston, B., Taylor, P., Stansby, P.: Solitary wave
  transformation, breaking and run-up at a beach.
\newblock In: {Proceedings of the Institution of Civil Engineers-Maritime
  Engineering}, vol. 159, pp. 97--105. Thomas Telford Ltd (2006)

\bibitem{BMBBF04}
Briganti, R., Musumeci, R.E., Bellotti, G., Brocchini, M., Foti, E.:
  {Boussinesq modeling of breaking waves: Description of turbulence.}
\newblock Journal of Geophysical Research: Oceans \textbf{109}(C07015) (2004)

\bibitem{Brun2018}
Brun, M.K., Kalisch, H.: {Convective wave breaking in the KdV equation}.
\newblock Analysis and Mathematical Physics \textbf{8}(1), 57--75 (2018)

\bibitem{hansen79}
Buhr~Hansen, J., Svendsen, I.: Regular waves in shoaling water, experimental
  data.
\newblock Tech. Rep.~21, Institute of Hydrodynamics and Hydraulic Engineering,
  Technical University of Denmark (1979)

\bibitem{Castro2005}
Castro, M., Ferreiro, A., Garc{\'\i}a-Rodr{\'\i}guez, J., Gonz{\'a}lez-Vida,
  J., Mac{\'\i}as, J., Par{\'e}s, C., V{\'a}zquez-Cend{\'o}n, M.: {The
  numerical treatment of wet/dry fronts in shallow flows: application to
  one-layer and two-layer systems}.
\newblock Mathematical and Computer Modelling \textbf{42}(3-4), 419--439 (2005)

\bibitem{CEA20123317}
Cea, L., V{\'a}zquez-Cend{\'o}n, M.: {Unstructured finite volume discretisation
  of bed friction and convective flux in solute transport models linked to the
  shallow water equations}.
\newblock Journal of Computational Physics \textbf{231}(8), 3317 -- 3339 (2012)

\bibitem{chanson2004}
Chanson, H.: {Hydraulics of open channel} flow.
\newblock Elsevier (2004)

\bibitem{Cienfuegos2010}
Cienfuegos, R., Barth{\'e}lemy, E., Bonneton, P.: Wave-breaking model for
  {B}oussinesq-type equations including roller effects in the mass conservation
  equation.
\newblock Journal of Waterway, Port, Coastal, and Ocean Engineering
  \textbf{136}(1), 10--26 (2009)

\bibitem{Delis2008}
Delis, A., Kazolea, M., Kampanis, N.: A robust high-resolution finite volume
  scheme for the simulation of long waves over complex domains.
\newblock International Journal for Numerical Methods in Fluids \textbf{56}(4),
  419--452 (2008)

\bibitem{dingemans1997}
Dingemans, M.: {Water Wave Propagation Over Uneven Bottoms: Linear wave
  propagation}.
\newblock Advanced series on ocean engineering. World Scientific Publishing Co
  Pte Ltd (1997)

\bibitem{favre1935}
Favre, H.: {\'E}tude th{\'e}orique et exp{\'e}rimentale des ondes de
  translation dans les canaux d{\'e}couverts.
\newblock Publications du Laboratoire de recherches hydrauliques annex{\'e}
  {\'a} l'{\'E}cole polytechnique f{\'e}d{\'e}rale de Zurich. Dunod (1935)

\bibitem{fbcr15}
Filippini, A., Bellec, S., Colin, M., Ricchiuto, M.: On the nonlinear behavior
  of {B}oussinesq type models: {A}mplitude-velocity vs amplitude-flux forms.
\newblock {Coastal Engineering} \textbf{99}, 109--123 (2015)

\bibitem{filippini2016}
Filippini, A., Kazolea, M., Ricchiuto, M.: A flexible genuinely nonlinear
  approach for wave propagation, breaking and run-up.
\newblock Journal of Computational Physics \textbf{310}, 381--417 (2016)

\bibitem{Harten}
Harten, A., Hyman, J.: Self adjusting grid methods for one-dimensional
  hyperbolic conservation laws.
\newblock Journal of Computational Physics \textbf{50}, 235--269 (1983)

\bibitem{KazoleaPM}
Kazolea, M.: Personal communication.

\bibitem{Kazolea2013}
Kazolea, M., Delis, A.: A well-balanced shock-capturing hybrid finite
  volume--finite difference numerical scheme for extended 1{D} {B}oussinesq
  models.
\newblock Applied Numerical Mathematics \textbf{67}, 167--186 (2013)

\bibitem{Kazolea14}
Kazolea, M., Delis, A.I., Synolakis, C.E.: Numerical treatment of wave breaking
  on unstructured finite volume approximations for extended {B}oussinesq-type
  equations.
\newblock J.Comput.Phys. \textbf{271}, 281--305 (2014)

\bibitem{kazolea2018}
Kazolea, M., Ricchiuto, M.: On wave breaking for {Boussinesq}-type models.
\newblock Ocean Modelling \textbf{123}, 16--39 (2018)

\bibitem{Kennedy2000}
Kennedy, A., Chen, Q., Kirby, J., Dalrymple, R.: {Boussinesq modeling of wave
  transformation, breaking and run-up. {I}: 1{D}}.
\newblock Journal of Waterway, Port, Coastal and Ocean Engineering
  \textbf{126}(1), 39--47 (2000)

\bibitem{Kermani}
Kermani, M., Plett, E.: {Modified entropy correction formula for the Roe
  scheme} p.~83 (2001)

\bibitem{book_lannes}
Lannes, D.: {The Water Waves Problem: Mathematical Analysis and Asymptotics}.
\newblock Mathematical surveys and monographs. American Mathematical Society
  (2013)

\bibitem{LeVeque}
LeVeque, R.J.: {Finite Volume Methods for Hyperbolic Problems}.
\newblock Cambridge Texts in Applied Mathematics. Cambridge University Press
  (2002)

\bibitem{LH69}
Longuet-Higgins, M.S.: On wave breaking and the equilibrium spectrum of
  wind-generated waves.
\newblock Proceedings of the Royal Society of London. A. Mathematical and
  Physical Sciences \textbf{310}(1501), 151--159 (1969)

\bibitem{Madsen1991}
Madsen, P.A., Murray, R., S{\o}rensen, O.R.: {A new form of the Boussinesq
  equations with improved linear dispersion characteristics}.
\newblock Coastal engineering \textbf{15}(4), 371--388 (1991)

\bibitem{Madsen1992}
Madsen, P.A., S{\o}rensen, O.R.: {A new form of the Boussinesq equations with
  improved linear dispersion characteristics. Part 2. A slowly-varying
  bathymetry}.
\newblock Coastal engineering \textbf{18}(3-4), 183--204 (1992)

\bibitem{Melville1988}
Melville, W., Rapp, R.J.: The surface velocity field in steep and breaking
  waves.
\newblock Journal of Fluid Mechanics \textbf{189}, 1--22 (1988)

\bibitem{Okamoto2006}
Okamoto, T., Basco, D.R.: {The Relative Trough Froude Number for initiation of
  wave breaking: Theory, experiments and numerical model confirmation}.
\newblock Coastal Engineering \textbf{53}(8), 675--690 (2006)

\bibitem{Pelanti}
Pelanti, M., Quartapelle, L., Vigevano, L.: A review of entropy fixes as
  applied to {R}oe's linearization.
\newblock Teaching material of the Aerospace and Aeronautics Department of
  Politecnico di Milano  (2001)

\bibitem{Ricchiuto2009}
Ricchiuto, M., Bollermann, A.: Stabilized residual distribution for shallow
  water simulations.
\newblock Journal of Computational Physics \textbf{228}(4), 1071--1115 (2009)

\bibitem{rf14}
Ricchiuto, M., Filippini, A.: Upwind residual discretization of enhanced
  {B}oussinesq equations for wave propagation over complex bathymetries.
\newblock Journal of Computational Physics \textbf{271}, 306--341 (2014)

\bibitem{ROEBER20121}
Roeber, V., Cheung, K.F.: Boussinesq-type model for energetic breaking waves in
  fringing reef environments.
\newblock Coastal Engineering \textbf{70}, 1--20 (2012)

\bibitem{Schaffer1993}
Sch{\"a}ffer, H.A., Madsen, P.A., Deigaard, R.: A boussinesq model for waves
  breaking in shallow water.
\newblock Coastal engineering \textbf{20}(3-4), 185--202 (1993)

\bibitem{shi2012}
Shi, F., Kirby, J.T., Harris, J.C., Geiman, J.D., Grilli, S.T.: {A high-order
  adaptive time-stepping TVD solver for Boussinesq modeling of breaking waves
  and coastal inundation}.
\newblock Ocean Modelling \textbf{43}, 36--51 (2012)

\bibitem{Mingham2009}
Shiach, J.B., Mingham, C.G.: {A temporally second-order accurate Godunov-type
  scheme for solving the extended Boussinesq equations}.
\newblock Coastal Engineering \textbf{56}(1), 32--45 (2009)

\bibitem{SKOTNER1999905}
Skotner, C., Apelt, C.: {Application of a Boussinesq model for the computation
  of breaking waves: Part 1: Development and verification}.
\newblock Ocean engineering \textbf{26}(10), 905--925 (1999)

\bibitem{Sorensen2004}
S{\o}rensen, O.R., Sch{\"a}ffer, H.A., S{\o}rensen, L.S.: Boussinesq-type
  modelling using an unstructured finite element technique.
\newblock Coastal Engineering \textbf{50}(4), 181--198 (2004)

\bibitem{Synolakis1987}
Synolakis, C.E.: The runup of solitary waves.
\newblock Journal of Fluid Mechanics \textbf{185}, 523--545 (1987)

\bibitem{Tissier2012}
Tissier, M., Bonneton, P., Marche, F., Chazel, F., Lannes, D.: {A new approach
  to handle wave breaking in fully non-linear Boussinesq models}.
\newblock Coastal Engineering \textbf{67}, 54--66 (2012)

\bibitem{Tonelli2009}
Tonelli, M., Petti, M.: {Hybrid finite volume--finite difference scheme for 2DH
  improved Boussinesq equations}.
\newblock Coastal Engineering \textbf{56}(5-6), 609--620 (2009)

\bibitem{Tonelli2010}
Tonelli, M., Petti, M.: {Finite volume scheme for the solution of 2D extended
  Boussinesq equations in the surf zone}.
\newblock Ocean Engineering \textbf{37}(7), 567--582 (2010)

\bibitem{Tonelli2011}
Tonelli, M., Petti, M.: {Simulation of wave breaking over complex bathymetries
  by a Boussinesq model}.
\newblock Journal of Hydraulic Research \textbf{49}(4), 473--486 (2011)

\bibitem{Tonelli2012}
Tonelli, M., Petti, M.: Shock-capturing boussinesq model for irregular wave
  propagation.
\newblock Coastal Engineering \textbf{61}, 8--19 (2012)

\bibitem{treske1994}
Treske, A.: {Undular bores (Favre-waves) in open channels-experimental
  studies}.
\newblock Journal of Hydraulic Research \textbf{32}(3), 355--370 (1994)

\bibitem{vmf15}
Viviano, A., Musumeci, R.E., Foti, E.: {A nonlinear rotational, quasi-2DH,
  numerical model for spilling wave propagation}.
\newblock Applied Mathematical Modelling \textbf{39}(3-4), 1099--1118 (2015)

\bibitem{Wei95}
Wei, G., Kirby, J.T., Grilli, S.T., Subramanya, R.: {A fully nonlinear
  Boussinesq model for surface waves. Part 1. Highly nonlinear unsteady waves}.
\newblock Journal of Fluid Mechanics \textbf{294}, 71--92 (1995)

\bibitem{zelt1991run}
Zelt, J.: The run-up of nonbreaking and breaking solitary waves.
\newblock Coastal Engineering \textbf{15}(3), 205--246 (1991)

\end{thebibliography}
%

\end{document}